\documentclass[11pt, reqno]{amsart}
\usepackage{amsmath}
\usepackage{amssymb}
\usepackage{amsfonts}
\usepackage[usenames]{color}
\usepackage{version}
\usepackage{graphics}

 \newtheorem{theorem}{Theorem}[section]
 \newtheorem{corollary}[theorem]{Corollary}
 \newtheorem{lemma}[theorem]{Lemma}
 \newtheorem{proposition}[theorem]{Proposition}

\newtheorem{definition}[theorem]{Definition}
\newtheorem{remark}[theorem]{Remark}
\newtheorem{example}[theorem]{Example}
\newtheorem{fact*}{Fact}

\DeclareMathOperator\hol{Hol}
\DeclareMathOperator{\RE}{Re}

\DeclareMathOperator{\diag}{diag}
\DeclareMathOperator\rank{rank}
\DeclareMathOperator\re{Re}

\newcommand\half{\tfrac 12}
\newcommand\dd{\mathrm d}

\newcommand{\s}{\mathcal{S}_2}

\newcommand\e{\mathrm{e}}
\newcommand\h{\mathcal{H}}

\newcommand{\M}{\mathcal{M}}

\newcommand{\T}{\mathbb{T}}

\newcommand{\F}{\mathcal{F}}
\newcommand{\D}{\mathbb{D}}
\newcommand{\C}{\mathbb{C}}
\newcommand{\G}{\mathcal{G}}
\newcommand{\rhp}{\mathbb H}

\newcommand{\schur}{\mathcal{S}_2}
\newcommand{\R}{\mathbb{R}}

\newcommand{\tr}{\operatorname{tr}}

\newcommand{\ip}[2]{\left\langle #1, #2 \right\rangle}

\newcommand{\inv}{^{-1}}

\renewcommand\L{\mathcal{L}}

\newcommand{\ph}{\varphi}

\newcommand\ga{\gamma}

\newcommand\la{\lambda}
\newcommand\La{\Lambda}

\newcommand\beq{\begin{equation}}

\newcommand\eeq{\end{equation}}
\newcommand\df{\stackrel{\rm def}{=}}

\newcommand\nn{\nonumber}
\newcommand\bbm{\begin{bmatrix}}
\newcommand\ebm{\end{bmatrix}}
\newcommand\bpm{\begin{pmatrix}}
\newcommand\epm{\end{pmatrix}}

\numberwithin{equation}{section}

\begin{document}

\title[$\mu$-synthesis as a semidefinite program]{A case of mu-synthesis as a quadratic semidefinite program }

\author{Jim Agler, Zinaida A. Lykova and N. J. Young}
\date{9th December, 2012}

\begin{abstract} 
We analyse a special case of the robust stabilization problem under structured uncertainty.
We obtain a new criterion for the solvability of the spectral Nevanlinna-Pick problem, which is a special case of the $\mu$-synthesis problem of $H^\infty$ control in which $\mu$ is the spectral radius.
Given $n$ distinct points $\la_1,\dots,\la_n$ in the unit disc and $2\times 2$ nonscalar complex matrices $W_1,\dots,W_n$,  the problem is to determine whether there is an analytic $2\times 2$ matrix function $F$ on the disc such that $F(\la_j)=W_j$ for each $j$ and the supremum of the spectral radius of $F(\la)$ is less than $1$ for $\la$ in the disc.  The condition is that the minimum of a quadratic function of pairs of positive $3n$-square matrices subject to certain linear matrix inequalities in the data be attained and be zero.
\end{abstract}

\maketitle

\section{Introduction} \label{intro}

In this paper we study an optimisation problem that arises in the design of a stabilizing controller for a linear time-invariant system that is subject to structured uncertainty.  We show that, in a special case, the existence of a robustly stabilizing controller is equivalent to the condition that the minimum of a quadratic objective function of a matrix pair subject to a linear matrix inequality (LMI) be attained and be zero.

Robust control theory provides a rigorous framework for the formulation and analysis of specifications of control systems for plants that are subject to sundry types of uncertainty; see for example \cite{Do,DFT} or  \cite[Chapter 8]{DuPa}.  One of the tools of the theory is the {\em structured singular value} of an operator or matrix corresponding to a given uncertainty class (\cite{Do} or \cite[Definition 8.13]{DuPa}); this is a cost function that generalizes the operator norm.  It is  denoted by $\mu$, and leads to the ``$\mu$-synthesis problem", which is a problem of optimization  over a class of analytic matrix functions in a disc or half-plane.   Special cases of the $\mu$-synthesis problem are the Nehari and Nevanlinna-Pick problems, which have classical solutions, but in virtually no other case is there an analytic solution -- there are only approximate numerical methods that have neither guaranteed convergence nor error bounds.

The $\mu$-synthesis problem is an interpolation problem for analytic matrix functions.  It is a familiar fact that robust stabilization leads to interpolation problems (\cite{Fr87,DuPa}).  For a given nominal plant and uncertainty class, the set of all stable closed-loop transfer functions can be parametrised, resulting in a class of analytic matrix functions $F$ that are subject to interpolation conditions.  To maximize the  uncertainty region about the nominal plant that can be simultaneously stabilized one must minimize over $F$ the quantity $\sup_{\la} \mu (F(\la))$, where $\la$ varies over a disc or halfplane and $\mu$ is a cost function that encodes structural properties of the uncertainty set.  In general the interpolation conditions that $F$ satisfies are of the ``model matching" type (Proposition \ref{modelmatch} below), but in this paper we restrict ourself to the case that the values of $F$ are specified at finitely many points.

The computation of $\mu$ for an arbitrary block-structured uncertainty is known to be NP-hard \cite{BYDM,TO}.  However, one familiar and easily-computed instance of the structured singular value is the spectral radius $r(\cdot)$ of a matrix; it corresponds to a one-dimensional uncertainty class.  Accordingly the following is an instance of the $\mu$-synthesis problem.  We denote the unit disc of the complex plane $\C$ by $\D$.\\

\noindent {\bf Problem SNP}  {\em Given distinct points $\la_1, \dots, \la_n \in \D$ and target matrices $W_1, \dots, W_n$ of type $k \times k$ find an analytic $k \times k$-matrix-valued function $F$ such that 
\[
F(\la_j)= W_j \quad \mbox{  for } j=1,\dots,n,  \mbox{  and }
\]
\[
\sup_{\la\in\D} r(F(\la)) \quad \mbox{ is minimized.}
\]
}
This problem is called the {\em spectral Nevanlinna-Pick problem} and has attracted attention as a test problem for $\mu$-synthesis  (\cite{BFT90,BFT1} and other papers by these authors).  There are heuristic algorithms that calculate approximate solutions to $\mu$-synthesis problems (\cite[Section 9.3.3]{DuPa} or \cite{mathworks}).  These algorithms are in current industrial use, but they are nevertheless widely regarded as not fully satisfactory.  In the absence of an adequate analytic theory it is even difficult to test their outputs for closeness to optimality.

In this paper we prove a solvability criterion for the spectral Nevanlinna-Pick problem in the $2\times 2$ case (that is, $k=2$).  The main result is the following.  

\begin{theorem}\label{NPspectral-main}
Let $\la_1, \dots, \la_n$ be distinct points in $\D$ and let $W_1,\dots, W_n$ be  $2\times 2$ complex matrices, none of them a scalar multiple of the identity.  Let $s_j= \tr W_j, \ p_j = \det W_j$ for each $j$ and let $z_1,z_2, z_3$ be any three distinct points in $\D$.  The following  three  conditions are equivalent.
\begin{enumerate}
\item There exists an analytic $2\times 2$ matrix function $F$ in $\D$ such that 
\beq \label{interpMat-main}
F(\la_j) = W_j \quad \mbox{ for }\quad j=1,\dots,n
\eeq
and 
\beq \label{specCond-main}
r(F(\la)) \leq 1 \quad \mbox{ for all } \quad \la \in\D;
\eeq
\item    the semidefinite program 
\beq\label{object}
\min  \quad (\tr N)^2-\tr(N^2)
\eeq
subject to
\begin{align}\label{constraints}
N&= [N_{i\ell, jk}]_{i,j=1, \, \ell,k=1}^{n,3} \geq 0, \notag \\
 M&=[M_{i\ell,jk}]_{i,j=1,\, \ell,k=1}^{n,3}\geq0,\notag \\
\left[1-\overline{\left(\frac{2z_\ell p_i-s_i}{2-z_\ell s_i}\right)} \frac{2z_k p_j-s_j}{2-z_k s_j}\right] 
& \ge \bbm(1-\bar z_\ell z_k) N_{i\ell,jk}\ebm + \bbm(1- \bar \la_i \la_j) M_{i\ell,jk}\ebm
\end{align}
is feasible and the minimum \eqref{object} is attained and has value zero;
\item  the semidefinite program in {\rm (2)} admits a feasible pair $(N,M)$ such that $\rank N\leq1$.

\end{enumerate}
\end{theorem}
Here the inequality sign denotes the usual partial order on the space of $3n$-square Hermitian matrices.
This result is a part of  Theorems \ref{NPspectral-criterion} and \ref {NPspectral-criterion-two}   which are established in Section \ref{criterion}.

Neither (2) nor (3) provides a convex program to resolve the existence of the desired function $F$ in (1).  In (2) the feasible region is convex, but the objective function $(\tr N)^2-\tr(N^2)$, though non-negative and quadratic, is not concave, so that it may fail to attain a minimum at an extreme point of the feasible region and it may have many local minima.  Condition (3) is a simple reformulation of (2), since it is easy to show that $(\tr N)^2-\tr(N^2)=0$ if and only if $\rank N \leq 1$  (Lemma \ref{rank1crit} below).  In (3) there is no objective function, but the  rank constraint means that the feasible set is not convex.  Thus neither formulation is an LMI problem, as treated for example in \cite{boyd}.   We do not know whether (2) or (3) can be the basis of an efficient numerical procedure. 
There are, however,  many papers in the literature on the numerical solution of optimization and feasibility problems of types that include (2) and (3), and it may be that some of the methods proposed will provide (in this special case) an effective alternative to algorithms in the current literature such as \cite{BFT1,DuPa,mathworks}.   Rank-constrained LMIs have been studied in \cite{OHM2004,RC2002,GI1994,GS1996}, while there are many algorithms for the optimization of smooth functions over convex sets \cite{mathworks}.  In the programs (2) and (3) the feasible pair $(N,M)$ may be constrained to lie in a set having a known prior bound (Proposition \ref{boundsNM} below).

Of course relaxation of the rank constraint in (3) yields a {\em necessary} condition for the solvability of a $2\times 2$ spectral Nevanlinna-Pick problem in the form of the feasibility of a true LMI: see Corollary \ref{necessary}.

The paper is organized as follows.  Section \ref{robust} describes the robust stabilization problem, presents a concrete example and outlines the reduction of the problem to a model matching problem. Section \ref{symmetrized} describes the symmetrized bidisc $\Gamma$ and its magic functions. It also describes the reduction of a $2\times 2$ spectral Nevanlinna-Pick problem to an interpolation problem in
the space $\hol(\D,\Gamma)$ of analytic functions from $\D$ to $\Gamma$.    Section \ref{example} gives a worked  example which is more general than the one in Section \ref{robust}.  This example, though illustrative, is limited to the case of systems with only two right-half-plane poles, and so motivates the need to develop an alternative approach.    To this end Section \ref{duality} presents a duality between the space $\hol(\D,\Gamma)$  and a subset of the Schur class $\s$ of the bidisc.
In Sections \ref{schurBidisc} and \ref{criterion}
we use Hilbert space models for functions in $\s$ to obtain necessary and sufficient conditions for solvability of the interpolation problem in the space $\hol(\D,\Gamma)$.  Section \ref{realiz-Phi} presents another approach to the realization of the relevant functions in $\s$. In Section \ref{matricial} we give some matricial formulations of the solvability criteria for the $2\times 2$ spectral Nevanlinna-Pick problem.  In Section \ref{all}  the procedure for constructing interpolating functions in $\hol(\D,\Gamma)$ developed in Section \ref{criterion} is shown to be general: in principle it yields {\em all} possible interpolants.   Section \ref{implement} contains some remarks about the numerical implementation of our procedure.

The closed unit disc in $\C$ will be denoted by $\D^-$ and the unit circle by $\T$.  

The complex conjugate transpose of a matrix $A$ will be written $A^*$.  The symbol $I$ will denote an identity operator or an identity matrix, according to context.  For a matrix $A$ and a non-zero scalar $\la$ we shall sometimes write $A/\la$ as a synonym for $\la^{-1} A$. The right half plane $\{s:\re s > 0\}$ will be denoted by $\rhp$, and $RH^\infty$ will be the space of real-rational matricial functions that are analytic and bounded on $\rhp$; the type of the matrices will be understood from the context.

\section{Robust stabilization} \label{robust}
The theory of the robust stabilization of a plant subject to structured uncertainty is a particularly elegant chapter of $H^\infty$ control.  It was developed in the 1980s and 1990s (see \cite{Do, Do85} and many references in \cite{DuPa}); it is well described in \cite[Chapters 8 and 9]{DuPa}.   In this section we sketch the reduction of the robust stabilization problem to an optimization problem of `model matching' type.   Consider the system $\Sigma$  in Figure 1:
\vspace*{0.5cm}

\hspace*{1.5cm} \includegraphics{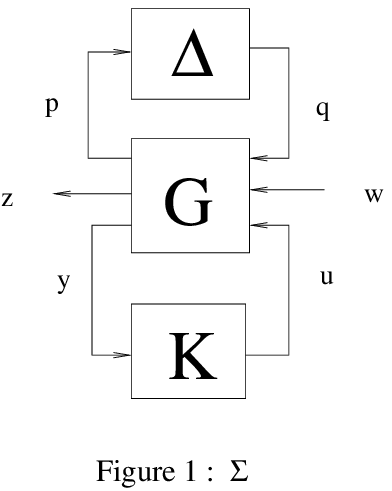} \hspace*{2cm}  \includegraphics{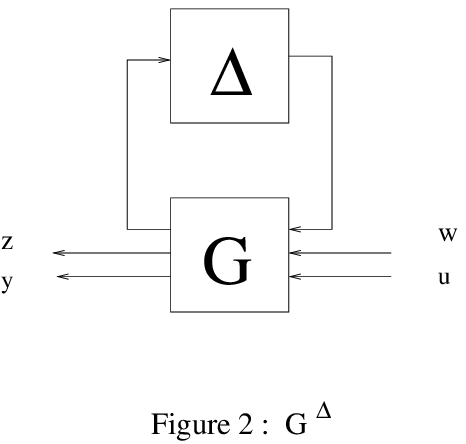} 

$\Sigma$ is a feedback system with uncertainty; $\Delta, G$ and $K$ denote finite-dimensional linear time-invariant systems, identified with their transfer functions, which are real rational matrix-valued functions of the frequency variable $s$. The nominal plant 
\beq\label{G}
G=\bbm G_{ij} \ebm_{i,j=1}^3
\eeq
is given and is proper.   Uncertainty is modelled by the assumption that the `true plant'  $G^\Delta$ is given by Figure 2 for some $\Delta$ which belongs to a prescribed `uncertainty set' $\mathbf \Delta$ but is otherwise unknown.  We shall say that $G$ is {\em robustly stabilizable with respect to} $\mathbf \Delta$ if there exists a stable controller $K$ such that the configuration in Figure 1 is stable for all $\Delta\in \mathbf \Delta$.  Mathematically the requirements are that 
\begin{enumerate}
\item the controller $K$ belong to $RH^\infty$, 
\item the system in Figure 1 be well posed for all $\Delta\in \mathbf \Delta$ and
\item   the system in Figure 1 be stable for all $\Delta\in \mathbf \Delta$.  
\end{enumerate}
A system is said to be {\em well posed} if  the transfer functions between different branches of the interconnection are well defined -- in particular, if all the inverses occurring in the transfer functions exist in the ring of square rational matrix functions of appropriate type.  See \cite[page 282]{DuPa} for a fuller discussion of well-posedness.

Here is an instance of the robust stabilization problem.  For $\ga\in\rhp$ let $b_\ga$ denote the stable allpass function of degree $1$ (or Blaschke factor)
\beq\label{defBlasch}
b_\ga(s)= \frac{s-\ga}{s+\bar\ga}, \quad s\in\rhp.
\eeq
Let $f, g \in RH^\infty$ be defined by
\beq\label{deffg}
f(s)=\frac{4}{3}\frac{s+3}{s+1}, \qquad g(s) = -\frac{(s+3)(s+5)}{3(s+1)^2}.
\eeq
It may be verified that
\beq\label{bezout}
fb_1+gb_3=1.
\eeq
As in \cite[page 262]{DuPa}, define the uncertainty set $\mathbf \Delta_{1,0}$ by
\beq\label{onedim}
\mathbf \Delta_{1,0} \df \{\delta I: |\delta | \leq 1\}.
\eeq
\begin{example}\label{firstex}
Let 
\beq\label{firstG}
G = \bbm G_{ij}\ebm_{i,j=1}^3=\bbm b_1b_3 & b_1+gb_3 & 1 & g & 1 & 0\\
	10b_3+fb_1&2b_{\sqrt{3}} + b_1b_3 & fb_1/b_3 & 1 & 0 & b_1/b_3 \\
	1& gb_3/b_1 & 1 & g/b_1 & 1/b_1 & 0 \\
	0 & 1& f/b_3&1&0&1/b_3\\
	1&0&1/b_3&0&0&b_1/b_3\\
	0&b_3/b_1&0&1/b_1& b_3/b_1&0 \ebm,
\eeq
regarded as partitioned into $2\times 2$ blocks.
Does there exist a robustly stabilizing controller $K$ for the system $\Sigma$ of Figure $1$ with respect to the uncertainty set $\mathbf \Delta_{1,0}$?
\end{example}
The answer to this question is given in Proposition \ref{answer} on page \pageref{answer}.

For a general plant $G$, evidently some stabilizability assumption on $G$ is a prerequisite for the existence of the $K$ we are seeking.
We shall  assume that $G$ is stabilizable, that is, there exists a controller $K$ such that the lower loop of Figure 1, which is the system $M(K)$  shown in Figure 3 below, is internally stable.  With this assumption, by \cite[Lemma 5.4]{DuPa}, $K$ stabilizes $G$ if and only if it stabilizes $G_{33}$.  However, stabilization of $G$ does not imply stabilization of $G^\Delta$ for a general $\Delta$.

In the notation of \cite[page 196]{DuPa},
\begin{align}\label{M}
M(K)= \bbm M_{ij} \ebm_{i,j=1}^2  &=\underline{S}(G,K) \notag\\
	&\df  \bbm G_{11} & G_{12} \\ G_{21} & G_{22} \ebm +\bbm G_{13}\\ G_{23} \ebm K(1-G_{33} K)\inv \bbm G_{31} & G_{32}\ebm.
\end{align}

\hspace*{1.7cm} \includegraphics{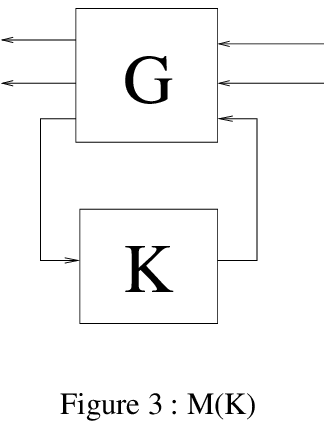} \hspace*{3cm} \includegraphics{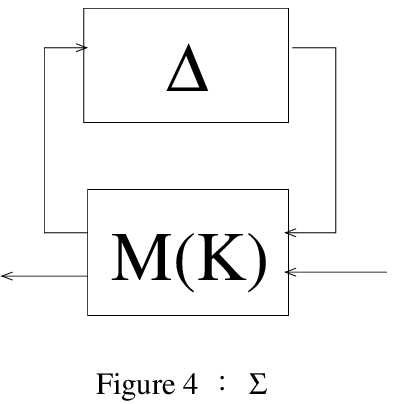}

The system $\Sigma$ of Figure 1 can also be written as in Figure 4, and so it has transfer function
\begin{align}\label{F}
 \overline{S}(M(K),\Delta) \df M_{22}+ M_{21}\Delta(1-M_{11}\Delta)\inv M_{12}, \quad \Delta\in\mathbf \Delta.
\end{align}
 For a fixed $K$, Theorems 8.22 and 9.8 of \cite{DuPa} assert that $ \overline{S}(M(K),\Delta )$ exists and belongs to $ RH^\infty$ for all $\Delta\in \mathbf \Delta$ if and only if $K$ stabilizes $G_{33}$ and
\beq\label{m11}
\sup_{s\in\rhp} \mu(M_{11}(s), \mathbf \Delta) < 1.
\eeq
Here $\mu(\cdot, \mathbf \Delta)$ denotes the structured singular value of a matrix relative to the uncertainty set $\mathbf \Delta$ \cite[Definition 8.13]{DuPa}.   For the present paper it is enough to assume that $\mathbf \Delta$ is such that $\mu(\cdot, \mathbf \Delta)=r(\cdot)$, the spectral radius.

Now the set of all stabilizing controllers of $G_{33}$ has an elegant parametrization due originally to Youla \cite[Theorems 5.13 and 5.14]{DuPa}.  Let $G_{33}$ have the doubly coprime factorization
\beq\label{doubly}
G_{33} = \hat N \hat M\inv= \tilde M\inv\tilde N
\eeq
 over $RH^\infty$ where  $\tilde N, \tilde M,\tilde X, \tilde Y, \hat N, \hat M,\hat X$ and $\hat Y$ belong to $RH^\infty$ and satisfy
\beq\label{bezout2}
\bbm \tilde X &- \tilde Y\\-\tilde N & \tilde M\ebm \bbm \hat M&\hat Y\\ \hat N& \hat X \ebm = I
\eeq
(every proper real rational plant admits such a factorization: see \cite[Proposition 5.10]{DuPa}).
Then the general stabilizing controller of $G_{33}$ is given by 
\beq\label{K}
K=(\hat Y- \hat M Q)(\hat X-\hat N Q)\inv = (\tilde X - Q\tilde N)\inv (\tilde Y-Q\tilde M)
\eeq
for some $Q\in RH^\infty$ such that $\hat X(\infty)- \hat N(\infty) Q(\infty)$ is invertible. 
 Moreover, for $K$ in equation \eqref{K}, 
\beq\label{mkq}
M(K)=\hat T_1-\hat T_2Q\hat T_3
\eeq
where
\begin{align}\label{that}
\hat T_1 &= \bbm G_{11} & G_{12} \\ G_{21} & G_{22} \ebm + \bbm G_{13}\\ G_{23} \ebm  \hat Y \tilde M \bbm G_{31} & G_{32} \ebm, \notag \\
\hat T_2 &= \bbm G_{13} \\ G_{23} \ebm \hat M,  \\
\hat T_3 &=  \tilde M \bbm G_{31} & G_{32} \ebm. \notag
\end{align}
 To summarize:
\begin{proposition}\label{modelmatch}
Let $G=\bbm G_{ij}\ebm_{i,j=1}^3$ be a stabilizable plant and let $G_{33}$ have the doubly coprime factorization \eqref{doubly}--\eqref{bezout2} over $RH^\infty$.  Let the zero matrix belong to $\mathbf{\Delta} \subset RH^\infty$.  There exists a controller $K\in RH^\infty$ such that the system $\Sigma$ of Figure $1$ is internally stable for all $\Delta\in\mathbf\Delta$ if and only if there exists $Q \in RH^\infty$ such that $\hat X(\infty)- \hat N(\infty) Q(\infty)$ is nonsingular and 
\beq\label{muQ}
\sup_{s\in\rhp} \mu\left((T_1-T_2QT_3)(s), \mathbf{\Delta}\right) < 1
\eeq
where
\beq\label{Tj}
T_1=G_{11}+G_{13} \hat Y \tilde M G_{31}, \quad T_2= G_{13}\hat M, \quad T_3= \tilde M G_{31}.
\eeq
Moreover the general robustly stabilizing controller of $\Sigma$ for the uncertainty set $\mathbf \Delta$ is given by equation  \eqref{K} for some $Q\in RH^\infty$ such that $\hat X(\infty)- \hat N(\infty) Q(\infty)$ is nonsingular and inequality \eqref{muQ} holds.
\end{proposition}
The problem of whether there exists a function $Q$ with these properties is called a {\em model matching problem}.
\begin{proof}
Suppose there does exist a controller $K$ that stabilizes $\Sigma$ for all $\Delta\in\mathbf\Delta$.  Since $0 I \in\mathbf \Delta$, in particular $K$ stabilizes the plant
\[
G^0= \bbm G_{22} & G_{23} \\ G_{32} & G_{33} \ebm.
\]
Hence $K$ stabilizes $G_{33}$, and so by the Youla parametrization there exists $Q\in RH^\infty$ such that
$\hat X(\infty)- \hat N(\infty) Q(\infty)$ is nonsingular and $K$ satisfies equations \eqref{K}.  Furthermore $K$ satisfies the inequality \eqref{m11}.  In view of equation \eqref{mkq},
\[
M_{11}= \bbm I&0\ebm M(K) \bbm I\\0\ebm = T_1-T_2 Q T_3
\]
where
\[
T_1= \bbm I & 0\ebm \hat T_1 \bbm I \\ 0\ebm, \qquad  T_2= \bbm I & 0\ebm\hat  T_2, \qquad T_3=\hat T_3 \bbm I \\ 0\ebm
\]
It follows from equations \eqref{that} that $T_1, T_2, T_3$ are given by equations \eqref{Tj}.  Thus necessity holds in Proposition \ref{modelmatch}.  To prove sufficiency one simply reverses the steps.
\end{proof}

By the example on \cite[page 257]{DuPa}  $\mu(\cdot, \mathbf \Delta_{1,0})=r(\cdot)$, the spectral radius.  Thus the robust stabilization problem for the nominal plant $G$ with uncertainty set $ \mathbf \Delta_{1,0}$ reduces to the following.

{\em Find $Q \in RH^\infty$ such that $\hat X(\infty)- \hat N(\infty) Q(\infty)$ is nonsingular and }
\beq\label{Qcond}
\sup_{s\in\rhp} r\left( (T_1-T_2QT_3)(s) \right) < 1.
\eeq
Now specialise further to the case that $T_2$ and $T_3$ are scalar matrix functions and the zeros of $T_2T_3$ in $\rhp$ are simple, say $s_1,\dots,s_n \in\rhp$.  Then the set of functions $T_1-T_2QT_3$ as $Q$ ranges through $RH^\infty$ is just the set of functions in $RH^\infty$ that agree with $T_1$ at $s_1,\dots,s_n$, and the question becomes whether there exists a function $F\in RH^\infty$ such that
\beq\label{rhpsnp}
F(s_j)=T_1(s_j) \quad  \mbox{ for } j=1,\dots,n\quad \mbox{ and }\quad \sup_{s\in \rhp} r(F(s)) < 1,
\eeq
which (after the application of a Cayley transform) is an instance of Problem SNP.  In the case that $G_{11}$ and $G_{33}$ are $2\times 2$ matrix functions Problem SNP can be analysed with the aid of the theory of the symmetrized bidisc, as presented in the next section.

\section{The symmetrized bidisc and its magic functions}\label{symmetrized}

The {\em open} and {\em closed symmetrized bidiscs} are the subsets
\beq\label{defG}
\G=\{(z+w,zw): |z| < 1, \, |w| < 1 \}
\eeq
and 
\beq\label{defGamma}
\Gamma = \{(z+w,zw): |z| \leq 1, \, |w| \leq 1 \}
\eeq
of $\C^2$.  They are relevant to the $2\times 2$ spectral Nevanlinna-Pick problem because, for a $2\times 2$ matrix $A$,
\[
r(A) < 1 \Leftrightarrow (\tr A, \det A) \in \G
\]
and
\beq\label{rGam}
r(A) \leq 1 \Leftrightarrow (\tr A, \det A) \in \Gamma.
\eeq
Accordingly, if $F$ is an analytic $2\times 2$ matrix function on $\D$ satisfying $r(F(\la)) \leq 1$ for all $\la \in\D$ then the function $(\tr F, \det F)$ belongs to the space $\hol(\D,\Gamma)$ of analytic functions from $\D$ to $\Gamma$.  A converse statement also holds: every $\ph\in\hol(\D,\Gamma)$ lifts to an analytic $2\times 2$ matrix function $F$ on $\D$ such that $(\tr F,\det F)= \ph$ and consequently $r(F(\la)) \leq 1$ for all $\la \in\D$ \cite[Theorem 1.1]{AY04T}.  The $2\times 2$ spectral Nevanlinna-Pick problem can therefore be reduced to an interpolation problem in $\hol(\D,\Gamma)$.  There is a slight complication in the case that any of the target matrices are scalar multiples of the identity matrix; for simplicity we shall exclude this case in the present paper. 

The relation \eqref{rGam} scales in an obvious way: for $\rho > 0$,
\[
r(A) \leq \rho \Leftrightarrow (\tr A, \det A) \in \rho \cdot\Gamma
\]
where 
\[
\rho\cdot(s,p) \df (\rho s,\rho^2 p) \quad \mbox{ and }\quad   \rho \cdot\Gamma\df \{\rho\cdot(s,p):(s,p)\in\Gamma\}.
\]

The following result is a refinement of \cite[Theorem 1.1]{AY04T}.
\begin{proposition}\label{reducetoGamma}
Let $\la_1,\dots,\la_n$ be distinct points in $\D$ and let $W_1,\dots,W_n$ be $2\times 2$ matrices, none of them a scalar multiple of the identity.  The following two statements are equivalent.
\begin{enumerate}
\item There exists a rational ${2\times 2}$ matrix function $F$, analytic in $\D$, such that 
\[
F(\la_j) = W_j \quad \mbox{ for } j=1,\dots, n
\]
and
\beq\label{rineq}
 \sup_{\la\in\D}r(F(\la)) < 1;
\eeq
\item there exists a rational function $h\in\hol(\D, \G)$ such that
\beq\label{hatnodes}
h(\la_j) = (\tr W_j, \det W_j) \quad \mbox{ for } j= 1,\dots,n,
\eeq
and $h(\D)$ is relatively compact in $\G$.
\end{enumerate}
\end{proposition}
\begin{proof}
(1)$\Rightarrow$(2)    Let $F$ be any function with the properties described in (1) and let $r_0$ be the supremum in the inequality \eqref{rineq}.  Let $f= (\tr,\det)\circ F$; then $f\in\hol(\D,\G)$, $f$ is rational and $f(\D) \subset r_0\cdot \Gamma\subset \G$, and so $f(\D)$ is relatively compact in $\G$.  

(2)$\Rightarrow$(1)   Let $h=(h_1,h_2)$ be as in (2).  Since $\G=\bigcup_{0<t<1} t\cdot\G$ and $h(\D)$ is relatively compact in $\G$ there exists $t\in (0,1)$ such that $h(\D) \subset t\cdot\G$.

Since the $2\times 2$ matrices $W_1,\dots,W_n$ are not scalar matrices, they are `nonderogatory', that is, their rational canonical forms have only one block, or alternatively, they are similar to companion matrices.  Hence there exist nonsingular matrices $P_1,\dots, P_n$ such that 
\[
 W_j=P_j\inv \bbm 0&1\\ -\det W_j & \tr W_j \ebm P_j, \qquad \mbox{ for } j=1,\dots,n.
\]
Pick a polynomial matrix $P(\la)$ such that $P(\la_j)=P_j$ for each $j$ and $P(\la)$ is nonsingular for every $\la\in\D$.  The matrix function
\[
F(\la)= P(\la)\inv \bbm 0&1\\ -h_2(\la)& h_1(\la) \ebm P(\la)
\]
is rational and analytic in $\D$.  Its characteristic polynomial is $z^2-h_1(\la)z+ h_2(\la)$ and so, since $h(\D) \subset t\cdot \Gamma$
\[
\sup_{\la\in\D} r(F(\la)) \leq t <1.
\]
\end{proof}

Certain simple rational functions play a central role in the analysis of $\Gamma$.

\begin{definition}\label{defPhi}
The function $\Phi$ is defined for $(z, s, p) \in \mathbb{C}^3$ such that $zs\neq 2$ by
\beq\label{defPhi-formula}
\Phi(z,s,p) = \frac{2zp-s}{2-zs} = -\half s+ \frac{(p-\tfrac 14 s^2)z}{1-\half s z}.
\eeq
\end{definition}

In particular, $\Phi$ is defined and analytic on $\mathbb{D} \times
\Gamma$ (since $|s| \le 2$ when $(s, p) \in \Gamma$),
 $\Phi$ extends analytically to  $(\Delta\times \Gamma) \setminus \{(z, 2\bar z, \bar z^2): z\in\T \}$.   See \cite{AY2} for an account of how $\Phi$ arises from operator-theoretic considerations.  The $1$-parameter family $\Phi(\omega, \cdot), \; \omega\in\T,$ comprises the set of {\em magic functions} of the domain $\G$.  The notion of magic functions of a domain is explained in \cite{AY08}, 
but for this paper all we shall need is the fact that
\[
\Phi(\D \times \Gamma) \subset \Delta
\]
and a converse statement: if $w \in \C^2$ and $|\Phi(z, w)| \le 1$ for all $z \in \D$ then $w \in \Gamma$;
see for example \cite[Theorem 2.1]{AY04} (the result is also contained in  \cite[Theorem 2.2]{AY1} in a different notation).

A {\em $\Gamma$-inner function} is the analogue for $\hol(\D,\Gamma)$ of inner functions in the Schur class.
A good understanding of rational $\Gamma$-inner functions is likely to play a part in any future solution of the finite interpolation problem for $\hol(\D, \Gamma)$, since such a problem has a solution if and only if it has a rational $\Gamma$-inner solution  (for example, \cite[Theorem 4.2]{Cost05} or Theorem \ref{NPspectral-criterion} below). 

\begin{definition}\label{Gam-in-funct} A  {\em $\Gamma$-inner function} is an analytic function $h : \D \to \Gamma$ such that the radial limit
\begin{equation}\label{radial}
\lim_{r \to 1-} h(r \lambda) \in b \Gamma
\end{equation}
for almost all $\lambda \in \T$, where $b\Gamma$ denotes the distinguished boundary of $\Gamma$.
\end{definition}
By Fatou's Theorem, the radial limit (\ref{radial}) exists for almost all 
 $\lambda \in \T$ with respect to Lebesgue measure. 
The  distinguished boundary $b\Gamma$ of $\G$ (or $\Gamma$)  is the \v{S}ilov boundary of the algebra of continuous functions on $\Gamma$ that are analytic in $\G$.  It is the symmetrisation of the 2-torus:
$$
b\Gamma= \{ (z+w,zw): |z|=|w|=1\}.
$$

\section{An example} \label{example}

Here is an example of a robust stabilization problem in which all the signals $p,q,z,w,y$ and $u$ in Figure 1 are two-dimensional.   The theory in  Section \ref{robust} reduces the robust stabilizability of this example to a two-point spectral Nevanlinna-Pick problem, which is amenable to a precise analysis.

\begin{example}\label {numEx}
Let $G= \bbm G_{ij}\ebm_{i,j=1}^3$ be a stabilizable plant such that, for some $a\in\C$ and $c\neq 0$, 
\begin{align}\label{formsG}
G_{11} &= \bbm b_1b_3 & b_1+gb_3 \\  ab_3+fb_1 & cb_{\sqrt{3}}+b_1b_3 \ebm, \nn\\
G_{13}&=\diag \{1, b_1/b_3\},\nn\\
G_{31}&=\diag\{ 1,b_3/b_1\}, \nn\\
G_{33} &= \bbm 0& b_1/b_3 \\ b_3/b_1 & 0 \ebm
\end{align}
where $f,g,b_\ga$ are as in equations \eqref{defBlasch} and \eqref{deffg}.

Find the values of $a,c$ for which there exists a robustly stabilizing controller for the system $\Sigma$ of Figure $1$ with respect to the uncertainty set 
$\mathbf \Delta_{1,0}$ of equation \eqref{onedim},
and describe the set of all robustly stabilizing controllers.
\end{example}

$G$ is unstable, having poles at $1$ and $3$ in the right halfplane.  There do exist stabilizable plants $G$ whose corner blocks are as in equations \eqref{formsG}; examples will be given later.  It will also transpire that the plant $G$ of Example \ref{firstex} is one such $G$ with $a=10, c=2$.
\begin{proposition}\label{solution}
There exists a robustly stabilizing controller in  Example {\rm \ref{numEx}} if and only if
\beq\label{bstrict}
|c| < \frac{1}{4-2\sqrt{3}}.
\eeq
\end{proposition}

\begin{proof}
Begin with a doubly coprime factorization of $G_{33}$.  The functions
\begin{align}\label{dbly}
\hat N&= \tilde N = \bbm 0&b_1\\b_3&0\ebm, \quad  \hat M=\diag\{b_1,b_3\},\quad \tilde M =\diag\{b_3,b_1\}, \notag \\
\tilde X &= \bbm f&-b_1 \\ -b_3 & g \ebm,   \qquad\tilde Y =  - \bbm b_3&g \\ f& b_1 \ebm, \\
\hat X&= \bbm g&-b_1\\ -b_3 & f \ebm,  \quad \qquad \hat Y= -\bbm  b_1&g \\ f& b_3 \ebm \notag
\end{align}
belong to $RH^\infty$ and satisfy equations \eqref{doubly} and \eqref{bezout2}.
In the notation of Proposition \ref{modelmatch}
\begin{align*}
T_1 &= G_{11}+G_{13} \hat Y \tilde M G_{31} 
	= \bbm 0 & b_1 \\ ab_3 & cb_{\sqrt{3}} \ebm, \\
T_2 &= G_{13}\hat M = b_1 I\\
T_3&=  \tilde M G_{31} = b_3 I.
\end{align*}

Proposition \ref{modelmatch} now asserts that $\Sigma$ is robustly stabilizable with respect to the uncertainty set $\mathbf \Delta_{1,0}$ if and only if there exists $Q\in RH^\infty$ such that 
\beq\label{Qinf}
\hat X(\infty)- \hat N(\infty) Q(\infty) \mbox{ is nonsingular and }
\eeq
\beq\label{supr}
\sup_{s\in\rhp} r\left( (T_1-b_1b_3 Q)(s) \right) < 1.
\eeq

Let $\kappa$ denote the Cayley transform
\[
\kappa: \D \to \rhp: \la \mapsto \frac{1+\la}{1-\la}
\]
The zeros $1,3$ of the scalar functions $T_2, T_3 \in\rhp$ correspond under $\kappa\inv$ to the points $0,\half$ respectively in $\D$. 

\begin{lemma}\label{zeta0}
Let $\la_1,\la_2,\zeta\in\D$.  There exists a rational function $h\in\hol(\D,\G)$ such that
\beq\label{ant}
h(\la_1)=(\zeta,0), \quad h(\la_2)=(-\zeta,0) \quad \mbox{ and } h(\D) \mbox{ is relatively compact in }\G
\eeq
if and only if
\[
|\zeta| < d(\la_1,\la_2)
\]
where $d$ denotes the pseudohyperbolic distance on $\D$.
\end{lemma}
\begin{proof}
Let $\delta_\G, C_\G$ be the Lempert function and Carath\'eodory distance on $\G$ respectively.  By definition of the Lempert function, for $z_1,z_2\in\G$ there exists $h\in\hol(\D,\G)$ such that $h(\la_1)=z_1$ and $h(\la_2)=z_2$ if and only if $\delta_\G(z_1,z_2) \leq d(\la_1,\la_2)$.  Moreover, if there is such an $h$, then in the case of $\G$ there is a {\em rational} $h$ with the same properties.  By \cite[Corollary 5.7]{AY04}, $\delta_\G= C_\G$, while \cite[Corollary 3.5]{AY04} gives an explicit formula for $C_\G$.  In particular, for any $\zeta\in\D$,
\beq\label{valdelta}
\delta_\G((\zeta,0), (-\zeta,0)) = C_\G((\zeta,0), (-\zeta,0)) = |\zeta|.
\eeq
$\Rightarrow$  Let $h\in\hol(\D,\G)$ satisfy conditions \eqref{ant}.  Then $h(\D) \subset t\cdot \G$ for some $t\in(0,1)$.  The function $g=t\inv\cdot h$ belongs to $\hol(\D,\G)$ and satisfies
\[
g(\la_1)= (t\inv\zeta,0), \quad g(\la_2)=(-t\inv\zeta,0).
\]
By definition of $\delta_\G$,
\[
\delta_\G((t\inv\zeta,0), (-t\inv\zeta,0))\leq d(\la_1,\la_2).
\]
Hence, by equation \eqref{valdelta},
\[
|t\inv \zeta| \leq d(\la_1,\la_2),
\]
and since $t<1$ it follows that $|\zeta| <  d(\la_1,\la_2)$.

$\Leftarrow$  Suppose that $|\zeta| <  d(\la_1,\la_2)$.
For $t\in (0,1]$,
\[
\frac{\dd}{\dd t} d(t\la_1,t\la_2)^2 = \frac{2t|\la_1-\la_2|^2}{|1-t^2\bar\la_2\la_1|^2 }\re \frac{1+t^2\bar\la_1\la_2}{1-t^2\bar\la_1\la_2} > 0.
\]
Hence there exists $t<1$ such that 
\[
\delta_\G((\zeta,0), (-\zeta,0))  = |\zeta|\leq d(t\la_1,t\la_2).
\]
  Consequently there is a rational function $g\in\hol(\D,\G)$ such that $g(t\la_1)=(\zeta,0), \ g(t\la_2) = (-\zeta,0)$.  Now the function $h(\la)=g(t\la)$ is rational, is analytic from $t\inv\D$ to $\G$ and maps $\la_1,\la_2$ to $(\pm\zeta,0)$.  Moreover $h(\D)\subset h(\D^-)$, a compact subset of $\G$.
\end{proof}
Suppose that $G$ is robustly stabilizable; then we may pick $Q\in RH^\infty$ such that conditions \eqref{Qinf} and \eqref{supr} hold.  Let $F=(T_1-T_2QT_3)\circ \kappa$. This $F$ is a rational analytic matrix function on $\D$ such that 
\begin{align}\label{Dsnp}
F(0)&= T_1(1)=\bbm 0&0\\ -\tfrac 12 a & (\sqrt{3}-2)c \ebm, \\
 F(\half) &= T_1(3)= \bbm 0&\half \\0& (2-\sqrt{3})c \ebm 
\end{align}
and
\beq\label{spec}
\sup_{\la\in\D} r(F(\la)) < 1.
\eeq
The rational function $h=(\tr,\det)\circ F$  belongs to $\hol(\D,\G)$,  $h(\D)$ is relatively compact in $\G$  (see Proposition \ref{reducetoGamma})  and
\begin{align}\label{2ptG}
h(0)&= (\tr F(0),\det F(0)) = ((\sqrt{3}-2)c,0), \notag\\
  h(\half) &= (\tr F(\half), \det F(\half)) = ((2-\sqrt{3})c, 0).
\end{align}
Hence by Lemma \ref{zeta0}, 
\begin{align*}
|(\sqrt{3}-2)c| &< d(0,\half)=\half,
\end{align*}
and therefore 
\beq\label{bcond}
|c| < \frac{1}{4-2\sqrt{3}} \approx 1.866.
\eeq

The condition \eqref{bcond} is also sufficient for the existence of a robustly stabilizing controller.  Suppose it is satisfied.  By Lemma \ref{zeta0} there exists  a rational function $h=(h_1,h_2) \in \hol(\D, C)$ for some compact subset $C$ of $\G$ such that the interpolation conditions \eqref{2ptG} hold.    Let
\[
h(1)=(\zeta,\eta) \in \Gamma.
\]
Since the matrices $T_1(1), T_1(3)$ in equations \eqref{Dsnp} are not scalar, there are nonsingular matrices $P_0,P_1$ such that  $P_0  T_1(1) P_0\inv, P_1 T_1(3)P_1\inv $ are companion matrices, that is,
\[
T_1(1) = P_0\inv \bbm 0&1\\ 0& (\sqrt{3}-2)c \ebm P_0, \quad 
T_1(3)= P_1\inv \bbm 0&1\\ 0& (2-\sqrt{3})c \ebm P_1.
\]
Let $V$ be a nonsingular matrix (to be chosen later) and let $P$ be a matrix polynomial such that $P(\la)$ is nonsingular for all $\la \in\D^-$ and
\[
P(0)=P_0, \quad P(\half)=P_1 \quad\mbox{ and } \quad P(1) = V.
\]
Define $F$ by
\[
F= P\inv \bbm 0&1\\ -h_2& h_1 \ebm P.
\]
Then
\beq\label{Fvals}
F(0) = T_1(1), \quad F(\half) = T_1(3) \quad \mbox{ and }\quad F(1) = V\inv \bbm 0&1\\ -\eta & \zeta \ebm V.
\eeq
Now let $Q\in RH^\infty$ be such that $F= (T_1-b_1b_3Q)\circ \kappa$.  On letting $\la\to 1$ (and hence $s=\kappa(\la) \to\infty$) in this relation we obtain
\[
V\inv \bbm 0&1\\ -\eta & \zeta \ebm V=F(1)=(T_1-b_1b_3Q)(\infty) = \bbm 0&1\\a&c\ebm-Q(\infty)
\]
and so
\beq\label{Qinfty}
Q(\infty) = \bbm 0&1\\a&c\ebm - F(1)= \bbm 0&1\\a&c\ebm - V\inv \bbm 0&1\\ -\eta & \zeta \ebm V.
\eeq
From the equations \eqref{dbly} we have
\begin{align*}
\hat X(\infty)&=  \bbm g&-b_1\\ -b_3 & f \ebm (\infty) = \bbm -\tfrac 13 & -1\\ -1& \tfrac 43 \ebm \\
\hat N(\infty)&= \bbm 0&b_1\\b_3&0\ebm (\infty) = \bbm 0&1\\1&0 \ebm.
\end{align*}
Hence
 \begin{align*}
\hat X(\infty) - \hat N(\infty) Q(\infty) &= \bbm -\tfrac 13 & -1\\- 1& \tfrac 43 \ebm - \bbm 0&1\\1&0 \ebm \left( \bbm 0&1\\a&c\ebm - V\inv \bbm 0&1\\ -\eta & \zeta \ebm V \right) \\
	&= \bbm 0&1\\1&0\ebm \left( V\inv \bbm 0&1\\ -\eta & \zeta \ebm V + Z\right)
\end{align*}
where
\[
Z=\bbm -1&\tfrac 13 \\ -\tfrac 13 -a& -1-c \ebm.
\]
Thus  $\hat X(\infty) - \hat N(\infty) Q(\infty)$  is nonsingular provided that
\beq\label{nonsing}
 \bbm 0&1\\ -\eta & \zeta \ebm  + VZ V\inv \mbox{ is nonsingular. }
\eeq
$Z$ is not a scalar matrix,  hence it is similar to its companion matrix
\[
\bbm 0&1\\ -(\frac {10}9+c+\tfrac 13 a) &- 2-c \ebm.
\]
Replacement of $Z$ by its companion form in equation \eqref{nonsing} shows that there is a $Q\in RH^\infty$ such that $\hat X(\infty) - \hat N(\infty) Q(\infty)$  is nonsingular if $V$ can be found such that 
\beq\label{mustbe}
\bbm 0&1\\-\eta&\zeta \ebm + V \bbm 0&1\\ -(\frac {10}9+c+\tfrac 13 a) & -2-c \ebm V\inv
\eeq
is nonsingular.  A suitable $V$ can be constructed whatever the values of $\zeta,\eta, a$ and $c$, and so $Q\in RH^\infty$ with the required properties exists.  On substituting this $Q$ into the formula \eqref{K} for $K$ we obtain the desired robustly stabilizing controller for $\Sigma$.

\end{proof}
To demonstrate that Example \ref{numEx} is not vacuous we need to show that the plant $G$ is stabilizable for a suitable choice of the second block row and column of $G$.  Choose $2\times 2$ matrix functions $R_{12}, R_{21}, R_{22}, R_{23}, R_{32} \in RH^\infty$ arbitrarily and let
\begin{align*}
G_{12} &=R_{12}+\bbm 0&g\\ fb_1/b_3 & 0 \ebm R_{32} \\
G_{21} &= R_{21}+ R_{23} \bbm 0&gb_3/b_1 \\0&0\ebm\\
G_{22} &= R_{22}+ R_{23} \bbm 0&g/b_1 \\f/b_3 &0\ebm R_{32}\\
G_{23}&=R_{23}\bbm b_1&0\\0&b_3 \ebm\inv\\
G_{32}&= \bbm b_3&0\\0&b_1\ebm\inv R_{32}.
\end{align*}
$G$ is now fully specified.   If we take $a=10, c=2$ and all the $R_{ij}$ equal to $I$ then we obtain the plant $G$ of Example \ref{firstex}.

To show that $G$ is stabilizable write down a right-coprime factorization of $G$ over $RH^\infty$.  It can be checked that $G=\hat n \hat m\inv$ where $\hat n=\bbm \hat n_{ij}\ebm$, 
\begin{align*}
\hat n_{11} &= \bbm b_1b_3 & b_1 \\ ab_3+fb_1 & cb_{\sqrt{3}}+b_1b_3 \ebm, \quad \hat n_{12} = R_{12}, \quad \hat n_{13}= b_1I, \\
\hat n_{21} &= R_{21}, \quad \hat n_{22} = R_{22}, \quad \hat n_{23} = R_{23}, \\
\hat n_{31} &= \bbm 1&0\\0&fb_3 \ebm, \quad \hat n_{32} = \bbm g&0\\0&f\ebm R_{32}, \quad \hat n_{33} =  \hat N = \bbm 0&b_1\\b_3&0\ebm,
\end{align*}
\[
\hat m= \bbm I & 0 & 0\\0&I&0\\ E_1&E_2 & \hat M \ebm
\]
and
\begin{align*}
E_1&=\bbm 0&-gb_3\\0&0\ebm, \qquad E_2 = -\bbm 0&g\\f&0\ebm R_{32}.
\end{align*}
Clearly $\hat n, \hat m \in RH^\infty$.  Moreover $\hat n, \hat m$ are right coprime, since
\[
\tilde x \hat m - \tilde y \hat n = I
\]
where 
\[
\tilde x= \bbm I&0&0\\0&I&0\\\tilde x_{31}& -R_{32} & \tilde X \ebm, \qquad \tilde y= \bbm 0&0&0\\0&0&0\\0&0&\tilde Y \ebm
\]
and
\[
\tilde x_{31}= - \bbm b_3&0\\ f&b_3 \ebm.
\]

We claim that the controller $K=\hat Y\hat X\inv$ stabilizes $G$ (and hence $G$ is stabilizable).     Since the functions defined in equations \eqref{dbly} satisfy the equation \eqref{bezout2} we have $\tilde M \hat X-\tilde N \hat Y=I$ and so $\hat X, \hat  Y$ are right coprime. According to \cite[Section 4.2, Theorem 1]{Fr87}, $K$ stabilizes $G$ if and only if the function
\[
\bbm \hat m & \bbm 0\\I\ebm \hat Y \\ \bbm 0&I\ebm \hat n & \hat X \ebm
\]
is invertible in $RH^\infty$ or equivalently, since $\bbm 0&I\ebm \hat n= \bbm G_{31}&G_{32}&G_{33}\ebm \hat m$, if and only if
\beq\label{8by8}
\bbm \hat m\inv&0\\0&I\ebm \bbm I&0&0&0 \\0&I&0&0\\0&0&I& \hat Y \\  G_{31}&G_{32}&G_{33}&\hat X \ebm\inv \in RH^\infty.
\eeq
It may be verified that
\[
\bbm I&0&0&0 \\0&I&0&0\\0&0&I& \hat Y \\  G_{31}&G_{32}&G_{33}&\hat X \ebm\inv  =\bbm I&0&0&0\\0&I&0&0\\ \hat Y C G_{31}&\hat Y C G_{32} & I+\hat YCG_{33}&-\hat Y C \\ 
	-CG_{31}& -CG_{32}& -CG_{33} & C \ebm
\]
where
\[
C= (\hat X-G_{33}\hat Y)\inv= \diag\{b_3,b_1\}.
\]
Hence $K$ stabilizes $G$ if and only if  
\[
\bbm I&0&0&0\\0&I&0&0\\ -\hat M\inv E_1 &-\hat M\inv E_2 &\hat M\inv &0\\0&0&0&I \ebm
   \bbm I&0&0&0\\0&I&0&0\\ \hat Y C G_{31}&\hat Y C G_{32} & I+\hat YCG_{33}&-\hat Y C \\ 
	-CG_{31}& -CG_{32}& -CG_{33} & C \ebm \in RH^\infty.
\]
It is a matter of straightforward calculation to check that this matrix product is
\[
\bbm I&0&0&0\\
	0&I&0&0\\
	-\bbm b_3&0\\f&b_3\ebm & -R_{32} & \tilde X & -\tilde Y\\
	-b_3 I & -R_{32}& -\tilde N & \tilde M \ebm
\]  
which belongs to $RH^\infty$.
Hence $G$ is stabilizable. 

The solution to Example \ref{firstex} now follows: it is the special case in which $a=10, c=2$ and all the $R_{ij}=I$.  Since Example \ref{numEx} admits a robustly stabilizing controller if and only if the inequality \eqref{bcond} holds, it is {\em not} robustly stabilizable when $c=2$. 
\begin{proposition}\label{answer}
The answer to the question posed in Example {\rm \ref{firstex}} is no. There is no robustly stabilizing controller for Example {\rm \ref{firstex}} with uncertainty set $\mathbf \Delta_{1,0}$.
\end{proposition}

It is striking that in Example \ref{numEx} the set of robustly stabilizing controllers does not depend on the five functions $R_{ij}$.  The proof of Proposition \ref{solution} contains the following.
\begin{proposition}\label{allcontrollers}
The robustly stabilizing controllers $K$ for the system $\Sigma$ of Example {\rm \ref{numEx}} when $  |c| < 1/(4-2\sqrt{3}) $ are the functions of the form
\[
K=(\hat Y- \hat M Q)(\hat X-\hat N Q)\inv = (\tilde X - Q\tilde N)\inv (\tilde Y-Q\tilde M)
\]
where the functions $\hat Y, \hat M, \hat X, \hat N, \tilde X, \tilde N, \tilde Y$ and $\tilde M$ are given by equations \eqref{dbly},
\[
Q=\left(\bbm 0&b_1\\ab_3 &cb_{\sqrt{3}}\ebm -F\circ \kappa\inv\right) /b_1b_3
\]
and $F$ is any rational analytic $2\times 2$ matrix function on $\D$ satisfying the conditions \eqref{Dsnp} and \eqref{spec}  and such that 
\[
F(1)+\bbm -1&\tfrac 13 \\ -\tfrac 13 -a&-1-c \ebm \quad \mbox{ is nonsingular.}
\]
\end{proposition}
Note that $F\circ\kappa\inv$ is proper for any choice of rational $F$.

Because Example \ref{numEx} is constructed so that the resulting model matching problem has scalar $T_2$ and $T_3$ and only two interpolation nodes, the function theory of $\Gamma$ is adequate for a full analysis of the robust stabilization problem.  For the general robust stabilization problem, even in the case of the uncertainty set $\mathbf \Delta_{1,0}$, currently known theory of $\Gamma$ does not suffice to decide the solvability of the corresponding model matching problem. 
The remainder of the paper develops an alternative approach that leads to the criterion in terms of a quadratic semidefinite program given in Theorem \ref{NPspectral-main}.

\section{Duality between $\hol(\D,\Gamma)$ and $\schur$}
\label{duality}
The Schur class of the bidisc will be denoted by $\s$:
\[
\s \df \hol(\D^2,\D^-).
\]

A strategy for the $2\times 2$ spectral Nevanlinna-Pick problem is as follows.
\begin{enumerate}
\item Reduce to an interpolation problem in $\hol(\D,\Gamma)$ as in Proposition \ref{reducetoGamma}.
\item The magic functions $\Phi(z, \cdot)$ induce a duality between $\hol(\D,\Gamma)$ and a subset of $\s$.
\item Use Hilbert space models for $\s$ to obtain a necessary and sufficient condition for solvability.
\end{enumerate}

For the second step, observe that since $\Phi(\D\times\Gamma) \subset \D^-$, if $h=(s,p)\in\hol(\D,\Gamma)$ then the function
\[
(z,\la) \mapsto \Phi(z,h(\la)) = \frac{2zp(\la)-s(\la)}{2-zs(\la)} \quad  \mbox{ for } z,\la\in\D
\]
belongs to $\s$. 
The simple observation that $\Phi$ induces a correspondence between $\hol(\D,\Gamma)$ and a subset of $\schur$ underlies the results in this paper. The correspondence was first developed in \cite{AYY}, where  a  realization theorem for $\hol(\D,\Gamma)$ was  proved and  some examples were calculated. In this section we answer the question:
\beq \label{Qu}
\text{ {\em which} subset of $\schur$ corresponds to $\hol(\D,\Gamma)$?}
\eeq
If $h=(s,p) \in \hol(\D,\Gamma)$ then, for any fixed $\la \in \D$, the map 
\beq 
z \mapsto \Phi(z, h(\la))= \frac{2zp(\la)-s(\la)}{2-zs(\la)}= \frac{2p(\la)z-s(\la)}{-zs(\la)+2}
\eeq
is a linear fractional self-map $f(z)= \frac{a z +b}{c z +d}$ of $\D$ with the property ``$b =c$". To make the last phrase precise, say that a linear fractional map $f$ of the complex plane has the property ``$b =c$" if $f(0)\neq \infty$ and either $f$ is a constant map or, for some $a, b$ and $d$ in $\C$,
\[
f(z)= \frac{a z +b}{b z +d}\; \; \text{for all} \;\; z \in \C \cup \{ \infty\}.
\]

The following is an easy calculation.

\begin{proposition} \label{propb=c} If $f$ is a non-constant linear fractional transformation then $f$ has the property ``$b =c$" if and only if $f(0)\neq \infty$ and
\[
f^{-1}(z) = - \frac{1}{f(-1/z)}\; \; \text{for all} \;\; z \in \C \cup \{ \infty\}.
\]
\end{proposition}

Here is an answer to the question \eqref{Qu}.

\begin{proposition} \label{Phi-lin-frac} 
Let $G$ be an analytic function on $\D^2$. There exists a function 
$h \in \hol(\D,\Gamma)$ such that 
\beq \label{GPhi}
G(z,\la) =\Phi(z, h(\la))\; \; \text{for all} \;\; z, \la \in \D 
\eeq
if and only if $G \in \schur$ and, for every $\la \in \D$,
$G(\cdot, \la)$ is a linear fractional transformation with the property ``$b =c$".
\end{proposition}
\begin{proof} 
Necessity is immediate (note that the relation \eqref{GPhi}
implies that $G(0, \la)= -\frac{1}{2} s(\la) \neq \infty$).

Conversely, suppose that $G \in \schur$ and $G(\cdot, \la)$ is a linear fractional transformation with the property ``$b =c$" for every $\la \in \D$.  Since $G(\cdot,\la)\in\schur$ for every $\la\in\D$, $G(\cdot,\la)$ does not have a pole at $0$, and therefore, by the ``$b=c$" property, we may write
\[
G(z,\la) = \frac{a(\la)z+b(\la)}{b(\la)z+1}=
b(\la) + \frac{(a(\la)-b(\la)^2)z}{b(\la)z+1} \quad \text{ for all } \;z,\la\in\D
\]
for some functions $a,b$ on $\D$.  Observe that this statement remains true even when $G(\cdot,\la)$ is constant for some $\la$.  Since
\[
b(\la) = G(0,\la)
\]
and
\[
a(\la)z= (b(\la)z + 1)\left(G(z,\la)-b(\la)\right) + b(\la)^2z,
\]
the functions $a$ and $b$ are analytic on $\D$. Let
\[
s(\la)= -2b(\la),\; p(\la)=a(\la) \quad \text{ for all } \;\la\in\D,
\]
and let $h =(s,p)$. Then $h$ is analytic on $\D$ and 
\[
\Phi(z, h(\la))=\frac{2zp(\la)-s(\la)}{2-zs(\la)}= \frac{a(\la)z+b(\la)}{b(\la)z +1}= G(z,\la). 
\]
Since $G \in \schur$ we have, for any $\la \in \D$
\[
\left|\Phi(z, h(\la))\right| \le 1 \quad \text{ for all } \;z \in\D,
\]
from which it follows that $h(\la) \in \Gamma$, see
\cite[Theorem 2.1 and Corollary 2.2]{AY04} and \cite[Proposition 3.2]{ALY11}. 
 Thus $h \in \hol(\D,\Gamma)$ has the required properties.
\end{proof}

\section{The Schur class of the bidisc}\label{schurBidisc}

Every function in $\s$ has a Hilbert space model \cite{Ag1}.  
That is to say, if $\ph\in\s$ then there exist a separable Hilbert space $\M$, a Hermitian projection $P$ on $\M$ and an analytic map $u:\D^2\to \M$ such that
\beq\label{defModel}
1-\overline{\ph(\mu)}\ph(\la) = \ip{(I-\mu_P^*\la_P)u(\la)}{u(\mu)} \quad \mbox{ for all } \la,\mu\in\D^2,
\eeq
where, for $\la=(\la_1, \la_2) \in \D^2$, $\la_P$ denotes $\la_1 P + \la_2 (I-P)$. 
This statement is contained in the proof of \cite[Theorem 1.12]{Ag1} -- see in particular equation (3.11).
 The triple $(\M,P,u)$ is called a {\em model} of $\ph$.

 The function $\Phi(z,h(\la))$ has the property that it is linear fractional in $z$ for every $\la\in\D$.  A consequence of this property is that the projection $P$ (corresponding to the variable $z$ in the defining equation \eqref{defModel}) has rank one for some model of the function, as we now show.

Denote by $\mathcal{S}^{2\times 2}$ the $2\times 2$ Schur class of the disc, that is, the set of analytic $2\times 2$ matrix functions $F$ on $\D$ such that $\|F(\la)\| \leq 1$ for all $\la\in\D$.  The following result (except for the uniqueness statement) is essentially \cite[Theorem 1.3]{AY04T}. 
\begin{proposition} \label{realiz-Gamma}
Let $h\in\hol(\D,\Gamma)$.  There exists  a unique function $F=\bbm F_{ij} \ebm \in \mathcal{S}^{2\times 2}$ such that 
\beq\label{trdet}
(\tr F, \det F) = h
\eeq and
\beq\label{propF}
F_{11}=F_{22}, \quad |F_{12}|=|F_{21}| \mbox{ a.e. on } \T, \quad F_{12} \mbox{ is either $0$ or outer and } F_{12}(0) \geq 0.
\eeq
Moreover, for all $\mu, \la \in \D$ and all $w, z \in \C$ such that 
$$
1 -   F_{22}(\mu)w \neq 0\;\;\text{and} \;\;
1 -  F_{22}(\la)z \neq 0,
$$ 
$F$ satisfies the identity  
\begin{eqnarray}
\label{1-PhiPhi-Pr}
1 -  \overline{\Phi(w, h(\mu))}\Phi(z, h(\la))
&=& (1-\bar{w}z)\overline{\gamma(\mu,w)}\gamma(\la,z) \nonumber\\ &~&\hspace{-2cm} 
+ \quad\eta(\mu, w)^* \left( I - F(\mu)^* F(\la)\right)\eta(\la, z),
\end{eqnarray}
where 
\begin{eqnarray}
\label{def-gamma-eta} 
\gamma(\la,z)&=&(1 - F_{22}(\la)z)^{-1} F_{12}(\la) \;\; \text{and} \nonumber\\
\eta(\la, z)&=& \bbm z\gamma(\la,z)\\ 1 \ebm.
\end{eqnarray} 
\end{proposition}

\begin{proof}
Let $h=(s,p)$.  Consider first the case that $s^2=4p$: then the function $F=\diag\{ \half s, \half s\}$ has the required properties \eqref{trdet} and \eqref{propF}, and moreover it is the only function with these properties.

In the case that $s^2 \neq 4p$ the $H^\infty$ function $\tfrac 14s^2 -p$ is nonzero and so it has an inner-outer factorization, expressible in the form
\[
\tfrac 14 s^2 - p = \ph \e^G
\]
where $\ph$ is inner, $\e^G$ is outer and $\e^G(0) \geq 0$.  Let
\beq\label{defF}
F=\bbm F_{ij} \ebm = \bbm \half s & \e^{\half G} \\ \ph \e^{\half G} & \half s \ebm.
\eeq
Then $\tr F=s$ and
\[
\det F = \tfrac 14 s^2- \ph\e^G = \tfrac 14 s^2 -(\tfrac 14 s^2 -p)=p.
\]
Clearly
\[
|F_{12}|= \e^{\re \half G} = |F_{21}| \quad\mbox{ a.e. on } \T,
\]
$F_{12}$ is outer and $F_{12}(0) > 0$.  Thus $F$ has the properties \eqref{trdet} and \eqref{propF}, and again it is easy to see that $F$ is the only function with these properties.

  We must show that $\|F\|_\infty \le 1$ on $\D$.  Let $f_1= F_{12}, \, f_2= F_{21}$.
At almost every point of $\T$
\begin{eqnarray}
\label{1-FF}
I -  F^*F
&=& I -\bbm \half \bar{s}  & \bar{f_2} \\ \bar{f_1} & \half \bar{s}  \ebm  \bbm \half s  & f_1 \\ f_2 & \half s \ebm \nonumber\\
\nonumber\\
&=& \frac{1}{4} \bbm 4 -|s|^2 - |s^2-4p| & -2\bar{s} f_1- 2 \bar{f_2} s\\ -2\bar{f_1}s- 2 f_2 \bar{ s} & 4 -|s|^2 - |s^2-4p| \ebm.
\end{eqnarray} 
The diagonal entries of the matrix on the right hand side are non-negative,
for if $s(\la) = z_1 + z_2 ,\; p(\la) = z_1 z_2$ where $z_1 , z_2 \in \D^- $, then
\begin{eqnarray}
\label{11entry}
 4 -|s|^2 - |s^2-4p|
&=& 4- |z_1 + z_2|^2 - |z_1 - z_2|^2 \nonumber\\
&=&  4- 2|z_1|^2 - 2|z_2|^2 \nonumber\\
&\ge& 0.
\end{eqnarray} 
Furthermore, for almost all $\la \in \T$,
\begin{eqnarray}
\label{12times21}
(-2\bar{s} f_1- 2 \bar{f_2} s)(-2\bar{f_1}s- 2 f_2 \bar{ s})
&=& |2\bar{s} f_1+ 2 \bar{f_2} s|^2\nonumber\\
&=& 4(|s|^2 |f_1|^2+ s^2 \bar{f_1} \bar{f_2} + \bar{s}^2 f_1 f_2+ |s|^2 |f_2|^2) \nonumber\\
&=& 4|s|^2 (|f_1|^2 +|f_2|^2) +8 \RE(\bar{s}^2 f_1 f_2) \nonumber\\
&=& 2|s|^2|s^2-4p| + 2\RE(\bar{s}^2 ( s^2 -4p)).
\end{eqnarray} 
Hence, for almost all  $\la \in \T$,
\begin{eqnarray}
\label{det(1-FF)}
16 \det(I -  F^*F)
&=& (4 -|s|^2 - |s^2-4p|)^2 - 2|s|^2|s^2-4p| - 2\RE(\bar{s}^2 (4p -s^2))\nonumber\\
&=& 16 + |s|^4 + |s^2 - 4p|^2 - 8|s|^2 - 8|s^2 - 4p|
+2|s|^2|s^2-4p|\nonumber\\
&~&\;\; - 2|s|^2|s^2-4p| - 2\RE(\bar{s}^2 ( s^2 -4p))\nonumber\\
&=& 16 - 8|s|^2 - 8|s^2 - 4p| \nonumber\\
&~&\;\; +|s|^4 + |s^2 - 4p|^2 
 - 2\RE(\bar{s}^2 ( s^2 -4p)).
\end{eqnarray} 
Note that
 $$
  |s|^4 + |s^2 - 4p|^2 - 2\RE(\bar{s}^2 (s^2 - 4p)) = |s^2 - (s^2 - 4p)|^2 = 16|p|^2 .
$$
Thus, for almost all $\la \in \T$,
$$
16 \det(I -  F^*F)= 16 + 16|p|^2 - 8|s|^2 - 8|s^2 - 4p|.
$$
Since $(s, p)$ maps $\D$ into $\Gamma$, by continuity $(s(\la), p(\la))$ can be written as $(z_1 +
z_2 , z_1 z_2 )$ for some $z_1 , z_2 \in \D^- $, and
\begin{eqnarray}
\label{det(1-FF)2}
16 \det(I -  F^*F)
&=& 16 + 16|z_1 z_2|^2 - 8|z_1 + z_2|^2 - 8|z_1 - z_2|^2\nonumber\\
&=& 16(1 + |z_1 z_2|^2 - |z_1|^2 -|z_2|^2)\nonumber\\
&=& 16 (1-|z_1|^2)(1-|z_2|^2)\ge 0.
\end{eqnarray} 

The inequalities (\ref{11entry}) and (\ref{det(1-FF)2}) show that
 $$      
I - F(\la)^* F(\la) \ge 0
$$
for almost all $\la \in \T$.  Thus $F\in\mathcal{S}^{2\times 2}$.

We now prove the identity \eqref{1-PhiPhi-Pr}.
For $\la \in \D$ and  for  $ z\in \C$ such that 
$1 -  \half s(\la)z \neq 0,$  
\begin{eqnarray}\label{proof-defPhi}
\Phi(z, h(\la)) 
&=&\frac{2zp(\la)-s(\la)}{2-zs(\la)}\nonumber\\
&=& -\half s(\la)+ \frac{(p(\la)-\tfrac 14 s(\la)^2)z}{1-\half s(\la) z}.
\end{eqnarray}

 Next apply a standard type of identity for linear fractional transformations; see \cite[Lemma 1.7]{AYY}. Let $H$, $U$ and $Y$ be Hilbert spaces. For any operator
\beq\label{defP}
P = \bbm P_{11}  & P_{12} \\ P_{21} &P_{22}  \ebm
: H \oplus U \to  H  \oplus Y,
\eeq
 denote by $\F_P$ the linear fractional transformation
\[
\F_P (X) = P_{22} + P_{21} X (I-  P_{11} X)^{-1}P_{12}
\]
defined for any operator $X$ on $H$ such that $I - P_{11}X$ is invertible. $\F_P (X)$, when
defined, is an operator from $U$ to $Y$. The following identity of standard type may
be verified by straightforward expansion.
Let 
$P = [P_{ij} ]^2_{i;j=1}$, $Q = [Q_{ij} ]^2_{i;j=1}$
be operators from $H \oplus U$ to $H \oplus Y $. For any pair of operators $X$, $Y$ on $H$ such that
$I - P_{11} X$ and $I - Q_{11} Y$ are invertible,
\begin{eqnarray}
\label{realization}
I - \F_Q (Y)^* \F_P (X) 
&=&   Q_{12}^* \left(I - Y^* Q^*_{11}\right)^{-1}
(I- Y^* X) \left(I - P_{11} X \right)^{-1} P_{12} \nonumber\\
&~& \hspace{-3cm}
+\bbm Q_{12}^* \left(I - Y^* Q^*_{11}\right)^{-1} Y^*  & I \ebm  
(I- Q^* P) 
\bbm X \left(I - P_{11} X \right)^{-1} P_{12}  \\ I \ebm .
\end{eqnarray} 
In the light of the definition \eqref{defF} of $F$,  for $\la \in \D$ and for any $ z \in \C$ such that 
$1 -  \half s(\la)z \neq 0,$
\begin{eqnarray}
\label{linear-trans-F}
\F_{F(\la)}(z) 
&=& \half s(\la) +f_1(\la) f_2(\la) z \frac{1}{1 -\half s(\la)z} \nonumber\\
&=& \half s(\la) + \frac{(\tfrac 14 s(\la)^2 - p(\la)) z}{1-\half s(\la) z},
\end{eqnarray}
and so, by equation (\ref{proof-defPhi}), 
\[
\F_{F(\la)}(z) = - \Phi(z, h(\la)).
\]
By the identity (\ref{realization}), for all $\mu, \la \in \D$ and for any $w, z \in \C$ such that 
$1 -  \half s(\mu)w \neq 0$ and  $1 -  \half s(\la)z \neq 0$,
 the expansion
\begin{eqnarray}
\label{1-PhiPhi}
1 -  \overline{\Phi(w, h(\mu))}\Phi(z, h(\la))
&=&  \overline{f_1(\mu)} (1 - \bar{w}\half \overline{s(\mu)})^{-1}(1-\bar{w}z)
(1 - \half s(\la)z)^{-1} f_1(\la)\nonumber\\
&~& \hspace{-2cm}
+\bbm \overline{f_1(\mu)} (1 - \bar{w}\half \overline{s(\mu)})^{-1} \bar{w} & 1 \ebm  
(I- F(\mu)^* F(\la)) 
\bbm z(1 - \half s(\la)z)^{-1} f_1(\la)  \\ 1 \ebm \nonumber\\
&=& (1-\bar{w}z)\overline{\gamma(\mu,w)}\gamma(\la,z) \nonumber\\
&~&\hspace{-2cm} +(1-\bar{\mu} \la) \eta(\mu, w)^* \frac{ I - 
F(\mu)^* F(\la)}{1- \bar{\mu} \la}\eta(\la, z)
\end{eqnarray}
holds, where $\gamma$ and $\eta$ are defined by equations \eqref{def-gamma-eta}.
\end{proof}

\section{An alternative proof of the
 realization of $\Phi(z, h(\la))$ on the bidisc}\label{realiz-Phi}

In this section we give an alternative proof of Proposition \ref{realiz-Gamma} in the special case that $h$ is $\Gamma$-inner.

\begin{proposition} \label{realiz-Gamma-inner} Let $h=(s,p):\D\to \Gamma$ be a 
 $\Gamma$-inner function. There exist a Hilbert space $\M$, an analytic function $F: \D \to \L(\C^2, \M)$ and an outer function $\gamma \in H^\infty$ such that 
\[ | \gamma (\la)|^2 = 1 -  |s(\la)|^2/4 \qquad \text{a. e. on} \; \T,
\]  
and for all $\mu, \la \in \D$ and for any $w, z \in \C$ such that 
$1 -  \half s(\mu)w \neq 0$ and $1 -  \half s(\la)z \neq 0,$ 
the identity
\begin{eqnarray}
\label{1-PhiPhi=PsiPsi-Pr}
1 - \overline{\Phi(w, h(\mu))} \Phi(z, h(\la))&=& (1-\bar{w} z) 
\left<\frac{ \gamma(\la)}{1 - \half  z s (\la)}, \frac{ \gamma(\mu)}{1 - \half  w s (\mu)} \right>
\nonumber\\
&~&  \hspace{-2cm} + (1 -\bar{\mu} \la) \left\langle F (\la) 
\begin{pmatrix}
1  \\ \\ \frac{z \gamma(\la)}{1 - \half z s (\la)}
\end{pmatrix}, F(\mu)
\begin{pmatrix}
1  \\ \\ \frac{w \gamma(\mu)}{1 - \half w s (\mu)}
\end{pmatrix} \right\rangle 
\end{eqnarray}
holds.
\end{proposition}

\begin{proof}
Since $(s,p)$ is $\Gamma$-inner, by \cite[Proposition 3.2(3)]{ALY11}, $s(\la)/p (\la)= \bar{s}(\la)$ and $|p(\la)|= 1$ for almost all $\la \in \T$. Thus, for almost all $\la \in \T$ and for $z \in \C$ such that 
$1 -  \half s(\la)z \neq 0,$ 
\begin{eqnarray}
\label{def-phi}
\Phi(z, h(\la)) 
&=& \frac{z - \half s(\la)/p(\la)}{1 - \half z s (\la)}p(\la) \nonumber\\
&=& \frac{z - \half \overline{s(\la)}}{1 - \half z s (\la)}p(\la). 
\end{eqnarray}
By the identity \eqref{def-phi}, since $|p(\la)|= 1$ for almost all $\la \in \T$,
 for  $w, z \in \C$ such that 
$
1 -  \half s(\mu)w \neq 0\;\;\text{and} \;\;
1 -  \half s(\la)z \neq 0,
$ 
 the expansion
\begin{eqnarray}
\label{1-phi-la} 
1-\overline{\Phi(w, h(\la))} \Phi(z,h(\la)) &=&
1- \overline{\frac{w - \half \overline{s(\la)}}{1 - \frac{1}{2} w s (\la)} p(\la)} \cdot \frac{z -\half \overline{s(\la)}}{1 - \frac{1}{2} z s (\la)}p(\la) \nonumber\\
\nonumber\\
&=&
1- \frac{\bar{w} - \half s(\la)}{1 - \half \bar{w} \overline{s(\la)}} \cdot \frac{z - \half \overline{s(\la)}}{1 - \half z s (\la)} \nonumber\\
 &=& \frac{\left(1 - |s(\la)|^2/4 \right) (1 -\bar{w} z)}
{\overline{\left(1 - \half s(\la) w \right)}\left(1 -  \half s(\la) z \right)} 
\end{eqnarray} 
holds  for almost all $\la\in \T$.
By a theorem of F. Riesz \cite{Du}, there exists an outer function $\gamma \in H^\infty$ such that 
\[
| \gamma (\la)|^2 = 1 -  |s(\la)|^2/4 \qquad \text{a. e. on} \; \T.
\]
Then, for  $w, z \in \C$ such that 
$
1 -  \half s(\mu)w \neq 0\;\;\text{and} \;\;
1 -  \half s(\la)z \neq 0,
$ 
\begin{eqnarray}
\label{1-phi-la-inner} 
\hspace{-1cm} 1- \overline{\Phi(w, h(\la))} \Phi(z,h(\la)) &=&
\left< \frac{\gamma(\la)}{1 - \half z s (\la)}, \frac{\gamma(\la)}{1 - \half w s (\la)} \right>_{\C} (1-\bar{w} z)\;\text{a. e. on}\; \T.
\end{eqnarray} 
Therefore, for all such $w, z \in \C$,
\begin{eqnarray}
\label{Gramians} 
\hspace{-1cm} 1+ \left< \frac{z\gamma(\la)}{1 - \half z s (\la)}, \frac{\bar{w} \gamma(\la)}{1 - \half w s (\la)} \right>_{\C} &=& 
\overline{\Phi(w, h(\la))} \Phi(z,h(\la)) +
 \nonumber\\
~ & ~ & \hspace{-1cm} \left< \frac{\gamma(\la)}{1 - \half z s (\la)}, \frac{\gamma(\la)}{1 - \half w s (\la)} \right>_{\C} \;\text{a. e. on}\; \T.
\end{eqnarray}
Thus the relation (\ref{Gramians}) can be expressed by the statement that,  for almost all $\la \in \T$ and for all $ z \in \C$ such that 
$1 -  \half s(\la)z \neq 0,$
 the Gramian  of the vectors 
\beq\label{gram3}
\begin{pmatrix}
1  \\ \\
\frac{z \gamma(\la)}{1 - \half z s (\la) }
\end{pmatrix} \in \C^2
\eeq
is equal to the Gramian of the vectors
\beq\label{gram4}
\begin{pmatrix}
\Phi(z, h(\la))  \\ \\ \frac{\gamma(\la)}{1 - \half z s (\la)}
\end{pmatrix} \in \C^2.
\eeq
Consequently there exists an isometric operator $L_{\la}$ on the vector space $\C^2$ which maps the vectors in equation \eqref{gram3} to the corresponding vectors in equation \eqref{gram4}. 

Define   $\Psi(\la)$, $\la \in \D$, by
\beq\label{defd-Psi}
\Psi(\la) \df \bbm -\half s (\la) & \frac{p(\la) - \frac{1}{4} s^2 (\la) }{\gamma(\la)} \\ \\  \gamma(\la) & \half s (\la) 
\ebm.
\eeq
Clearly $\Psi$ is analytic on $\D$. It is easy to check that, for all $\la \in \D$ and 
for any $ z \in \C$ such that $1 -  \half s(\la)z \neq 0$,
\begin{eqnarray}
\label{Psi-of-la}
\Psi(\la)
\begin{pmatrix}
1  \\ \\ \frac{z \gamma(\la)}{1 - \half z s (\la)}
\end{pmatrix} 
& = &
\bbm -\half s (\la) & \frac{p(\la) - \frac{1}{4} s^2 (\la) }{\gamma(\la)} \\ \\ \gamma(\la) & \half s (\la) \ebm
\begin{pmatrix}
1  \\ \\ \frac{z \gamma(\la)}{1 - \half z s (\la)} \end{pmatrix} 
\nonumber\\
& =& \begin{pmatrix}
\Phi(z, h(\la)) \\ \\ \frac{\gamma(\la)}{1 - \half z s (\la)}
\end{pmatrix}.
\end{eqnarray}
Since the boundary values $L_{\la}$ of $\Psi$ are isometries for almost all $\la \in \T$, the function $\Psi$
is in the Schur class $\mathcal{S}_{2 \times 2}$.

For all $\mu, \la \in \D$ and for any $w, z \in \C$ such that 
$1 -  \half s(\mu)w \neq 0\;\;\text{and} \;\;
1 -  \half s(\la)z \neq 0,$
\begin{eqnarray}
\label{1-PsiPsi}
\ip{\left(I - \Psi(\mu)^* \Psi(\la)\right)
\begin{pmatrix}
1  \\ \\ \frac{z \gamma(\la)}{1 - \half z s (\la)}
\end{pmatrix}} {
\begin{pmatrix}
1  \\ \\ \frac{w \gamma(\mu)}{1 - \half w s (\mu)}
\end{pmatrix} } 
&~& \nonumber\\
&~& \hspace{-9cm}
= \ip
{\begin{pmatrix}
1  \\ \\ \frac{z \gamma(\la)}{1 - \half z s (\la)}
\end{pmatrix}} {
\begin{pmatrix}
1  \\ \\ \frac{w \gamma(\mu)}{1 - \half w s (\mu)}
\end{pmatrix}}
-\ip{
\begin{pmatrix}
\Phi(z, h(\la))  \\ \\ \frac{ \gamma(\la)}{1 - \half z s (\la)}
\end{pmatrix}}{
\begin{pmatrix}
\Phi(w, h(\mu))  \\ \\ \frac{ \gamma(\mu)}{1 - \half w s (\mu)}
\end{pmatrix} }  \nonumber\\
&~& \hspace{-9cm}
= 1 - \overline{\Phi(w, h(\mu))} \Phi(z, h(\la)) - (1-\bar{w} z) 
\ip{\frac{ \gamma(\la)}{1 - \half  z s (\la)}}{ \frac{ \gamma(\mu)}{1 - \half  w s (\mu)} } .
\end{eqnarray}

Since  $\Psi$
is in the Schur class $\mathcal{S}_{2 \times 2}$, there exists a Hilbert space $\M$ and an analytic $F: \D \to \L(\C^2, \M)$ such that, for all $\mu, \la \in \D$, 
\[
I - \Psi(\mu)^* \Psi(\la)= (1 -\bar{\mu} \la) F(\mu)^* F (\la).
\]
The dimension of $\M$ is equal to the rank of $ \left[\frac{I - \Psi(\mu)^* \Psi(\la)}{1 -\bar{\mu} \la} \right]$.
Therefore, for all $\mu, \la \in \D$ and  for any $w, z \in \C$ such that 
$1 -  \half s(\mu)w \neq 0\;\;\text{and} \;\;
1 -  \half s(\la)z \neq 0,$ 
\begin{eqnarray}
\label{1-PhiPhi=PsiPsi}
1 - \overline{\Phi(w, h(\mu))} \Phi(z, h(\la))&=& (1-\bar{w} z) 
\ip{\frac{ \gamma(\la)}{1 - \half  z s (\la)}}{ \frac{ \gamma(\mu)}{1 - \half  w s (\mu)} }
\nonumber\\
&~& \hspace{-2cm}
+ \ip{\left(I - \Psi(\mu)^* \Psi(\la)\right)
\begin{pmatrix}
1  \\ \\ \frac{z \gamma(\la)}{1 - \half z s (\la)}
\end{pmatrix}}
{\begin{pmatrix}
1  \\ \\ \frac{w \gamma(\mu)}{1 - \half w s (\mu)}
\end{pmatrix} } 
\nonumber\\
&=&  (1-\bar{w} z) 
\left<\frac{ \gamma(\la)}{1 - \half  z s (\la)}, \frac{ \gamma(\mu)}{1 - \half  w s (\mu)} \right>
\nonumber\\
&~&  \hspace{-2cm} + (1 -\bar{\mu} \la) \left< F (\la) 
\begin{pmatrix}
1  \\ \\ \frac{z \gamma(\la)}{1 - \half z s (\la)}
\end{pmatrix}, F(\mu)
\begin{pmatrix}
1  \\ \\ \frac{w \gamma(\mu)}{1 - \half w s (\mu)}
\end{pmatrix} \right> .
\end{eqnarray}
\end{proof}

\begin{remark}\label{F-Psi} \rm  Let $h=(s,p):\D\to \Gamma$ be a  $\Gamma$-inner function. The relation between the analytic function $\Psi$ defined by \eqref{defd-Psi}
from Proposition \ref{realiz-Gamma-inner} and 
the analytic function $F$ defined by \eqref{defF}
from Proposition \ref{realiz-Gamma} is the following.
Recall that, for $\la \in \D$, 
\[
\Psi(\la) \df \bbm -\half s (\la) & \frac{p(\la) - \frac{1}{4} s^2 (\la) }{\gamma(\la)} \\ \\  \gamma(\la) & \half s (\la) 
\ebm
\]
where  $\gamma $ is  an outer function in  $H^\infty$ such that 
\[ 
| \gamma (\la)|^2 = 1 -  |s(\la)|^2/4 \qquad \text{a. e. on} \; \T.
\]
Since $(s,p)$ is $\Gamma$-inner, $s(\la)= \bar{s}(\la)p (\la)$ for almost all $\la \in \T$, by \cite[Proposition 3.2(3)]{ALY11}.
Therefore $\bar{s}^2(\la)p(\la) = |s(\la)|^2$ for almost all $\la \in \T$. Hence 
\[
| \gamma (\la)|^4= \left(1 -  |s(\la)|^2/4 \right)^2 = |p(\la) - s^2(\la)/4|^2 \qquad \text{a. e. on} \; \T.
\]
Thus, for
\[
f_1 = - \frac{p - \frac{1}{4} s^2  }{\gamma} \quad \text{and}\quad f_2 = \gamma,
\]
the functions $f_1 , f_2 \in  H^\infty$, $f_1 f_2 =  s^2 /4-p $ and
$$      
|f_1 (\la)| = |f_2 (\la)| \;\;    \text{ a. e. on}\;\; \T.
$$
Therefore
\[
F= \bbm \half s  & f_1 \\ f_2 &\half s  \ebm   = 
\bbm -1 & 0 \\ 0 & 1  \ebm
\bbm -\half s  & \frac{p - \frac{1}{4} s^2  }{\gamma} \\ \\  \gamma & \half s \ebm = \bbm -1 & 0 \\ 0 & 1  \ebm \Psi.
\]
\end{remark}

\section {Criteria for the solvability of spectral Nevanlinna-Pick problems}\label{criterion}

The following result (in combination with Proposition \ref{reducetoGamma}) contains Theorem \ref{NPspectral-main}. 
\begin{theorem}\label{NPspectral-criterion}
  Let $n\geq 1$, let $\la_1, \dots, \la_n$ be distinct points in $\D$  and let $(s_j,p_j)\in \Gamma$ for $j=1,\dots,n$.   Let $z_1,z_2,z_3$ be distinct points in $\D$.  The following  five  conditions are equivalent.
\begin{enumerate}
\item There exists an analytic function $h:\D\to \Gamma$ satisfying 
\beq\label{interp}
h(\la_j) = (s_j,p_j) \quad \mbox{ for }\quad  j=1,\dots,n;
\eeq
\item there exists a rational {\em $\Gamma$-inner} function $h$ satisfying equations \eqref{interp};
\item
there exist positive $3n$-square matrices  $N= [N_{i\ell,jk}]_{i,j=1,\, \ell,k=1}^{n,3}$ of rank at most $1$ and $M=[M_{i\ell,jk}]_{i,j=1,\, \ell,k=1}^{n,3}  $ such that,
for $1\le i,j \le n$ and $ 1 \le \ell,k \le 3$,
\beq\label{1.1}
1-\overline{\left(\frac{2z_\ell p_i-s_i}{2-z_\ell s_i}\right)} \frac{2z_k p_j-s_j}{2-z_k s_j}
 = (1-\bar z_\ell z_k) N_{i\ell,jk} + (1- \bar \la_i \la_j) M_{i\ell,jk};
\eeq
\item there exist positive $3n$-square matrices $N= [N_{i\ell, jk}]_{i,j=1, \, \ell,k=1}^{n,3}$ of rank  at most $1$ and $M=[M_{i\ell,jk}]_{i,j=1,\, \ell,k=1}^{n,3}  $ such that
\beq\label{inq1.1}
\left[1-\overline{\left(\frac{2z_\ell p_i-s_i}{2-z_\ell s_i}\right)} \frac{2z_k p_j-s_j}{2-z_k s_j}\right] \ge 
\bbm(1-\bar z_\ell z_k) N_{i\ell,jk}\ebm + 
\bbm(1- \bar \la_i \la_j) M_{i\ell,jk}\ebm;
\eeq
\item   the semidefinite program
\[
\min  \quad (\tr N)^2-\tr(N^2)
\]
subject to the linear matrix inequality \eqref{inq1.1} and the positivity conditions
\begin{align}\label{constraints2}
N&= [N_{i\ell, jk}]_{i,j=1, \, \ell,k=1}^{n,3} \geq 0, \notag \\
 M&=[M_{i\ell,jk}]_{i,j=1,\, \ell,k=1}^{n,3}\geq0,\notag 
\end{align}
   is feasible, attains its minimum and has value zero.

\end{enumerate}
\end{theorem}
\begin{proof} (2)$\Rightarrow$(1) and (3)$\Rightarrow$(4) are trivial. \\

(3)$\Rightarrow$(2)  Suppose $N,M$ as described exist.  Since $N$ is positive and of rank $1$ there exist scalars $\gamma_{jk}$ for $j= 1,\dots,n, \ k=1,2,3$ such that
\[
  N_{i\ell,jk} =  \bar \ga_{i\ell}\ga_{jk}.
\]
Likewise, since $M\geq 0$ there exist a Hilbert space $\M$ of dimension at most $3n$ and vectors $v_{jk} \in \M$ such that
\[
M_{i\ell,jk} = \ip{v_{jk}}{v_{i\ell}}_{\M}.
\]
Thus the relation \eqref{1.1} can be expressed by the statement that the Gramian of the vectors
\beq\label{gram1}
\begin{pmatrix}
-\Phi(z_k,s_j,p_j)  \\ \ga_{jk} \\ v_{jk}
\end{pmatrix} \in \C^2\oplus \M, \quad j=1,\dots,n, \quad k=1,2,3,
\eeq
is equal to the Gramian of the vectors
\beq\label{gram2}
\begin{pmatrix}
1  \\  z_k \ga_{jk} \\  \la_j v_{jk}
\end{pmatrix} \in \C^2\oplus \M, \quad j=1,\dots,n, \quad k=1,2,3.
\eeq
Consequently there exists a unitary operator $L$ on the finite-dimensional vector space $\C^2\oplus \M$ which maps the vectors in the expression \eqref{gram2} to the corresponding vectors in the expression \eqref{gram1}.  Write $L$ as a block operator matrix 
\[
L=\bbm A & B \\ C & D \ebm
\]
where $A, D$ act on $\C^2, \M$ respectively.  Then for $j=1,\dots,n, \ k=1,2,3$ we obtain the pair of equations
\begin{align*}
\bpm -\Phi(z_k,s_j,p_j)  \\ \ga_{jk} \epm &= A \bpm 1 \\ z_k\ga_{jk} \epm + B\la_j v_{jk}  \\
	v_{jk} &= C \bpm 1 \\ z_k\ga_{jk} \epm + D\la_j v_{jk}.
\end{align*}
From the second of these equations,
\beq\label{defvjk}
v_{jk}= (I-D\la_j)\inv C\bpm 1\\ z_k\ga_{jk} \epm,
\eeq
and hence
\beq\label{propABCD}
\bpm -\Phi(z_k,s_j,p_j)  \\ \ga_{jk} \epm = \left( A+B\la_j (I-D\la_j)\inv C\right) \bpm 1 \\ z_k\ga_{jk} \epm.
\eeq
Let
\beq\label{defPsi}
\Psi(\la)=  A+B\la (I-D\la)\inv C=\bbm a(\la) & b(\la) \\ c(\la) & d(\la) \ebm.
\eeq
Since $L$ is unitary and $\M$ is finite-dimensional, $\Psi$ is a rational $2 \times 2$  inner function, and hence 
\beq\label{defh}
h \df (\tr \Psi, \det \Psi)
\eeq
is a rational $\Gamma$-inner function.

We claim that $h$ satisfies the interpolation conditions \eqref{interp}.  By equation \eqref{propABCD},
\[
\bpm -\Phi(z_k,s_j,p_j)  \\ \ga_{jk} \epm = \Psi(\la_j) \bpm 1 \\ z_k\ga_{jk} \epm =
	\bpm a(\la_j) +b(\la_j)z_k \ga_{jk} \\ c(\la_j) + d(\la_j)z_k \ga_{jk} \epm
\]
for $j= 1,\dots,n, \ k=1,2,3$.
Eliminate $\ga_{jk}$ from these two equations to obtain
\[
\Phi(z_k,s_j,p_j) = -a(\la_j) - b(\la_j) z_k  (1-d(\la_j)z_k)\inv c(\la_j).
\]
That is to say that, for each $j\in\{1,\dots,n\}$, the linear fractional maps
\[
-\half s_j + \frac{(p_j - \tfrac 14 s_j^2)z}{1-\half s_j z} \quad \mbox{ and } \quad
	 -a(\la_j) - \frac{ b(\la_j) c(\la_j)z}{1-d(\la_j)z}
\]
agree at three distinct values of $z\in\D$.   It follows that the two maps are identical, which is to say that
\[
a(\la_j)= \half s_j, \quad b(\la_j)c(\la_j) = \tfrac 14 s_j^2 - p_j,  \quad  d(\la_j) = \half s_j.
\]
Hence
\[
\tr \Psi(\la_j) = a(\la_j) + d(\la_j) = s_j 
\]
 and 
\[
\quad \det \Psi(\la_j) = (ad-bc)(\la_j) = \tfrac 14 s_j^2 - (\tfrac 14 s_j^2 - p_j) = p_j.
\]
Thus $h(\la_j)= (s_j,p_j)$ for $j=1,\dots, n$ as required.  Thus (3)$\Rightarrow$(2).\\

(4)$\Rightarrow$(1) The proof of this statement is similar to that of (3) $\Rightarrow$(2). The only difference is that
 the relation \eqref{inq1.1} can be expressed by the statement that the Gramian of the vectors
\beq\label{gram1-inq}
\begin{pmatrix}
-\Phi(z_k,s_j,p_j)  \\ \ga_{jk} \\ v_{jk}
\end{pmatrix} \in \C^2\oplus \M, \quad j=1,\dots,n, \quad k=1,2,3,
\eeq
is {\em less than or equal} to the Gramian of the vectors
\beq\label{gram2-inq}
\begin{pmatrix}
1  \\  z_k \ga_{jk} \\  \la_j v_{jk}
\end{pmatrix} \in \C^2\oplus \M, \quad j=1,\dots,n, \quad k=1,2,3.
\eeq
Consequently there exists a {\em contractive} operator $L$ on the finite-dimensional vector space $\C^2\oplus \M$ which maps the vectors in the expression \eqref{gram2-inq} to the corresponding vectors in the expression \eqref{gram1-inq}.
Since $L$ is { contractive}, $\Psi$ defined by \eqref{defPsi} is in the {\em  $2 \times 2$ Schur class}, and hence 
\beq\label{defh-inq}
h = (\tr \Psi, \det \Psi) \in \hol (\D, \Gamma).
\eeq
The proof that $h(\la_j)= (s_j, p_j)$ is unchanged.

(1)$\Rightarrow$(3)  Suppose there exists an analytic function $h=(s,p):\D\to \Gamma$ satisfying equations \eqref{interp}. By Proposition \ref{realiz-Gamma},
there exists  an analytic function $F$ 
\[
F = \bbm \half s  & f_1 \\ f_2 &\half s  \ebm: \D \to \C^{2\times 2}
\]
such that $\|F \|_\infty \le  1$ on $\D$ and,
for all $\mu, \la \in \D$ and for any $w, z \in \C$ such that 
$1 -  \half s(\mu)w \neq 0$ and $1 -  \half s(\la)z \neq 0,$ 
\begin{eqnarray}
\label{1-PhiPhi-Th}
1 -  \overline{\Phi(w, h(\mu))}\Phi(z, h(\la))
&=& (1-\bar{w}z)\overline{\gamma(\mu,w)}\gamma(\la,z) \nonumber\\
&~&\hspace{-2cm} +(1-\bar{\mu} \la) \eta(\mu, w)^* \frac{ I - F(\mu)^* F(\la)}{1- \bar{\mu} \la}\eta(\la, z),
\end{eqnarray}
where 
\begin{eqnarray}
\label{def-gamma-eta-Th} 
\gamma(\la,z)&= &(1 - \half s(\la)z)^{-1} f_1(\la) \;\; \text{and} \nonumber\\
\eta(\la, z)& = & \bbm \gamma(\la,z)z\\ 1 \ebm.
\end{eqnarray} 
By assumption, for the given $\la_j \in \D$, $j =1, \dots, n$,
\beq\label{proof-interp}
h(\la_j) = (s_j,p_j) \quad \mbox{ for }\quad  j=1,\dots,n.
\eeq
Let $\mu=\la_i$ and $\la=\la_j$, $1 \le i,j \le n$,
in  (\ref{1-PhiPhi}). For all $w,z \in \D$, 
\begin{eqnarray}
\label{1-PhiPhi-ij}
1 -  \overline{\Phi(w, s_i, p_i)}\Phi(z, s_j, p_j)
&=&
1 -  \overline{\Phi(w, h(\la_i))}\Phi(z, h(\la_j))
\nonumber\\
&=& (1-\bar{w}z)\overline{\gamma(\la_i,w)}\gamma(\la_j,z) \nonumber\\
&~&\hspace{-2cm} +(1-\bar{\la_i} \la_j) \eta(\la_i, w)^* \frac{ I - F(\la_i)^* F(\la_j)}{1- \bar{\la_i} \la_j}\eta(\la_j, z).
\end{eqnarray}
Let $w =z_\ell, z =z_k$, $1 \le \ell, k \le 3$ in (\ref{1-PhiPhi-ij}). Since no $(s_j,p_j)$ is equal to $(2z_k, z_k^2)$ for any $k$, 
\begin{eqnarray}
\label{1-PhiPhi-ij-lk}
1 -  \overline{\Phi(z_\ell, s_i, p_i)}\Phi(z_k, s_j, p_j)
&=& (1-\bar{z_\ell}z_k)\overline{\gamma(\la_i,z_\ell)}\gamma(\la_j,z_k) \nonumber\\
&~&\hspace{-2cm} +(1-\bar{\la_i} \la_j) \eta(\la_i, z_\ell)^* \frac{ I - F(\la_i)^* F(\la_j)}{1- \bar{\la_i} \la_j}\eta(\la_j, z_k).
\end{eqnarray}
Since $\|F \|_\infty \le  1$ on $\D$, the matricial kernel
\[
(\la, \mu) \mapsto \frac{ I - F(\mu)^* F(\la)}{1- \bar{\mu} \la}: \D^2 \to M_2(\C)
\]
is positive.
Hence the positive $3n$-square matrices 
\[
N = [N_{i\ell,jk}]_{i,j=1,\, \ell,k=1}^{n,3} \df \left[\overline{\gamma(\la_i,z_\ell)}\gamma(\la_j,z_k)\right]_{i,j=1,\, \ell,k=1}^{n,3} 
\] 
of rank at most $1$ and 
\[
M=[M_{i\ell,jk}]_{i,j=1,\, \ell,k=1}^{n,3} \df \left[
\eta(\la_i, z_\ell)^* \frac{ I - F(\la_i)^* F(\la_j)}{1- \bar{\la_i} \la_j}\eta(\la_j, z_k)\right]_{i,j=1,\, \ell,k=1}^{n,3}  
\]
satisfy,
for $1\le i,j \le n$ and $ 1 \le \ell,k \le 3$,
\beq\label{proof1.1}
1-\overline{\Phi(z_\ell,s_i,p_i)}\Phi(z_k,s_j,p_j) = (1-\bar z_\ell z_k) N_{i\ell,jk} + (1- \bar \la_i \la_j) M_{i\ell,jk}.
\eeq

Therefore (1)$\Rightarrow$(3).

For the equivalence of (4) and (5) we require a simple observation involving exterior powers of matrices.  Recall that for any $r\times s$ matrix $A=\bbm a_{ij} \ebm$, the second exterior power $\bigwedge^2 A$ is an $\binom{r}{2} \times \binom{s}{2}$ matrix whose entries are the $2\times 2$ minors of $A$.  It follows that $A$ has rank  at most $1$ if and only if $\bigwedge^2 A =0$.
\begin{lemma}\label{rank1crit}
\begin{enumerate}
\item  For any self-adjoint matrix $A$,
\[
2\tr \bigwedge {}^2 A = (\tr A)^2- \tr(A^2).
\]
\item   If $A$ is a positive matrix then
\[
\rank A \leq 1 \quad \Leftrightarrow \quad (\tr A)^2- \tr(A^2) = 0 \quad \Leftrightarrow \quad (\tr A)^2- \tr(A^2) \leq 0.
\]
\end{enumerate}
\end{lemma}
\begin{proof} (1)
Let $A= \bbm a_{i,j} \ebm$.  Since $A=A^*$,
\begin{align*}
\tr (A^2) &= \tr (A^*A) = \sum_{i,j}|a_{ij}|^2 = \sum_i a_{ii}^2 + 2\sum_{i<j}|a_{ij}|^2,
\end{align*}
and so
\begin{align*}
2\tr \bigwedge {}^2 A &= 2\sum_{i<j} (a_{ii}a_{jj} - |a_{ij}|^2)\\
	&=2\sum_{i<j} a_{ii}a_{jj} - \left( \tr(A^2)- \sum_i a_{ii}^2\right) \\
	&= (\sum_i a_{ii})^2 - \tr(A^2) \\
	&= (\tr A)^2- \tr(A^2).
\end{align*}

(2)  Let $A\geq 0$.  Then also $\bigwedge^2 A \geq 0$ (if $A=B^*B$ then by the Cauchy-Binet formula $\bigwedge^2 A=\bigwedge^2 (B^*B)= (\bigwedge^2 B^*)(\bigwedge^2 B)=(\bigwedge^2 B)^*(\bigwedge^2 B) \geq 0$).
Thus
\begin{align*}
\rank A \leq 1 & \Leftrightarrow \bigwedge{}^2 A = 0 \\
	&\Leftrightarrow \tr \bigwedge{}^2 A \leq 0 \\
	& \Leftrightarrow  (\tr A)^2- \tr(A^2) \leq 0\\
	& \Leftrightarrow  (\tr A)^2- \tr(A^2) = 0.
\end{align*}
\end{proof}

(4)$\Leftrightarrow$(5)   The statement in (4) that there is a feasible pair $(N,M)$ for the LMI \eqref{inq1.1} with $N$ of rank $1$  means, in view of Lemma \ref{rank1crit}, that there is a feasible pair for which 
$(\tr N)^2- \tr(N^2) = 0$.  Since $(\tr N)^2- \tr(N^2) \geq 0$ whenever $N\geq 0$, it follows that the program in (5) is feasible and has value $0$.  The argument is reversible, and so (4)$\Leftrightarrow$(5).

Thus statements (1) to (5) are all equivalent.

\end{proof}

\begin{remark} \rm
(1) A natural choice of the points $z_1,z_2,z_3$ is $-1, 0, 1$.  This choice is not permitted in the theorem as stated above, which requires the $z_k$ to belong to $\D$.  However, the same proof works for $z_k\in\D^-$ provided that the denominators $2-z_ks_j$ in equation \eqref{1.1} are nonzero, which is so provided that no $(s_j,p_j)$ is equal to $(2z_k, z_k^2)$ for any $k$.\\

\noindent (2)  The matrix $\bbm 1- \overline{\Phi(z_\ell, s_i,p_i)}\Phi(z_k,s_j,p_j)\ebm$ in statements (3) and (4) of 
Theorem \ref{NPspectral-criterion} is positive if and only if the points $(s_j, p_j)$ are all equal and lie on the variety $s^2=4p$.

\noindent To see (2) study the proof of (3)$\Rightarrow$(2) with $\ga_{i\ell}=0, v_{i\ell}=0$.

 \noindent (3) The objective function $(\tr N)^2-\tr(N^2)$, though non-negative and quadratic, is not concave. For example, for $N=\lambda I \ge 0$ on $n$-dimensional space, since $\bigwedge^2 I$ is the identity on $\binom{n}{2}$-dimensional space,
\[
 \left({\tr} (\lambda I)\right)^2 - \tr\left((\lambda I)^2 \right)= 2 {\tr} \left(\lambda^2 I\right) = 2 \binom{n}{2} \lambda^2,
\]
which is not concave on $\R^+ I$.  The objective cannot be expected to attain its minimum at an extreme point of the feasible region.
\end{remark}

The next theorem relates the criterion of Theorem \ref{NPspectral-criterion} to the $\mu$-synthesis problem.  It also includes, for comparison, a statement of another criterion obtained previously in \cite{AY04T}.
\begin{theorem}\label{NPspectral-criterion-two}
Let $\la_1, \dots, \la_n$ be distinct points in $\D$ and let $W_1,\dots, W_n$ be  $2\times 2$ complex matrices, none of them a scalar multiple of the identity.  Let $s_j= \tr W_j, \ p_j = \det W_j$ for each $j$ and let $z_1,z_2, z_3$ be any three distinct points in $\D$.  The following four conditions are equivalent.
\begin{enumerate}
\item There exists an analytic $2\times 2$ matrix function $F$ in $\D$ such that 
\beq \label{interpMat}
F(\la_j) = W_j \quad \mbox{ for }\quad j=1,\dots,n
\eeq
and 
\beq \label{specCond}
r(F(\la)) \leq 1 \quad \mbox{ for all } \quad \la \in\D;
\eeq
\item there exists an analytic function $h : \D \to \Gamma$  such that 
\beq \label{tr-det}
h(\la_j ) = (\tr W_j ,
       \det W_j ),\;\; j = 1, 2, . . . , n;
\eeq
\item  there exists a bounded analytic $2\times 2$ matrix function $F$ in $\D$ such that conditions \eqref{interpMat} and \eqref{specCond} are satisfied, and in addition, both eigenvalues of $F(\la)$ have modulus $1$ for all $\la\in\T$;
\item there exist $b_1, \dots, b_n, c_1, \dots, c_n \in \C$ such that 
\beq
\label{criterion1}
\left[\frac{I - \bbm \half s_i  & b_i \\ c_i & -\half s_i  \ebm^*
\bbm \half s_j   & b_j \\ c_j & -\half s_j  \ebm}{1 -\bar{\la_i} \la_j} \right]^{n}_{i,j=1} \ge 0
\eeq
and
\beq
\label{crit1+cond2}
b_j c_j = p_j - \frac{s_j^2}{4}, \;\;  j=1,\dots, n.
\eeq
\end{enumerate}
\end{theorem}
\begin{proof} By \cite[Theorem 1.1 and Main Theorem 0.1]{AY04T}, conditions (1), (2) and (4) are equivalent. The equivalence (2)$\Leftrightarrow$(3) follows from (1)$\Leftrightarrow$(2) of Theorem \ref{NPspectral-criterion} and the fact that, for a $2 \times 2$ matrix $A$, both eigenvalues of $A$ have modulus $1$ if and only if $ (\tr A,  \det A)\in b \Gamma$.
\end{proof}

\begin{remark} \rm It is obvious that, under the hypothesis that none of the $W_j$ is a scalar multiple of the identity matrix, the conditions (1)-(4) of Theorem \ref{NPspectral-criterion-two} and the conditions (1)-(4)
of Theorem \ref{NPspectral-criterion} are equivalent.
\end{remark} 

\begin{remark} \rm In \cite[Main Theorem 0.1]{AY04T}, the authors proved that (1) and (4) of Theorem \ref{NPspectral-criterion} are equivalent under a genericity condition: none of the $W_j$ is a scalar multiple of the identity matrix. In \cite[Theorem 1.1]{Ber03} H. Bercovici removed this genericity condition and replaced condition (4) by a very similar one. 
\end{remark}

\section{Matricial formulations of the solvability criterion} \label{matricial}
There are more matricial ways of expressing the solvability criteria of Theorem \ref{NPspectral-criterion}.  Here are some of them.     Just for this section we shall denote by $I_n$ the $n\times n$ identity matrix. 
\begin{theorem}\label{NPspectral-LMI}
Let $\la_1, \dots, \la_n$ be distinct points in $\D$ and let $W_1,\dots, W_n$ be  $2\times 2$ complex matrices, none of them a scalar multiple of the identity.  Let $s_j= \tr W_j, \ p_j = \det W_j$ for each $j$.  Let $z_1,z_2$ and $z_3$ be distinct points of $\D^-$ such that no $(s_j,p_j)$ is equal to $(2z_k,z_k^2)$ for any $k$.

Let $3n$-square matrices $X, Z$ and $\Lambda$ be defined by
\begin{align}
X&= \left[ 1-\overline{\left(\frac{2z_\ell p_i - s_i}{2-z_\ell s_i}\right)}\frac{2z_kp_j -s_j}{2-z_k s_j}\right]_{i,j=1,\ell,k=1}^{n,3}, \label{defX}\\
\La &= \diag\{ \la_i\}_{i=1,\ell=1}^{n,3}, \label{defLa}\\
Z&= \diag\{z_\ell\}_{i=1,\ell=1}^{n,3}. \label{defZ}
\end{align}
The following  conditions are equivalent.
\begin{enumerate}
\item There exists an analytic $2\times 2$ matrix function $F$ in $\D$ such that 
\beq \label{interpMat-LMI}
F(\la_j) = W_j \quad \mbox{ for }\quad j=1,\dots,n
\eeq
and 
\beq \label{specCond-LMI}
r(F(\la)) \leq 1 \quad \mbox{ for all } \quad \la \in\D;
\eeq

\item  there exist positive $3n$-square matrices $N,M$ such that $\rank N \leq 1$ and
\beq\label{Xgeq}
X \geq N- Z^*NZ + M - \La^*M\La;
\eeq
\item the same as {\rm (2)} but for the replacement of $\geq$ by $=$;

\item there exist a positive $3n$-square matrix $M$, a $1\times 3n$ vector $\ga$ and a matrix $P$ of type $3n \times 2$
such that 
\beq\label{1.1bis-LMI}
\bbm -1 & 0 & 0 \\ 0 & 1 & 0\\ 0 & 0 & X \ebm \ge
 \bbm I_2 & 0 \\ P & I_{3n}\ebm
\bbm -1 & 0 & \ga \\ 0 & 1 & \ga Z\\  \ga^* &  Z^*\ga^* & M - \Lambda^* M  \Lambda \ebm
\bbm I_{2} & P^* \\ 0 & I_{3n} \ebm;
\eeq
\item  the same as {\rm (4)} but for the replacement of $\geq$ by $=$;
\item    the semidefinite quadratic program
\[
\min  \quad (\tr N)^2-\tr(N^2)
\]
subject to the conditions $N\geq 0, M \geq 0$ and the linear matrix inequality \eqref {Xgeq} is feasible and has value $0$. 
\end{enumerate}
\end{theorem}
Note that in $N,M,\La$ and $Z$ the rows are indexed by the pair $(i,\ell)$ and the columns by the pair $(j,k)$, where $i$ and $j$ run from $1$ to $n$, and $\ell$ and $k$ run from $1$ to $3$. 
\begin{proof}
The equivalences (1)$\Leftrightarrow$(2)$\Leftrightarrow$(3)$\Leftrightarrow$(6) are just reformulations of (1)$\Leftrightarrow$(3)$\Leftrightarrow$(4)$\Leftrightarrow$(5)  of Theorem \ref{NPspectral-criterion}. 

(2)$\Rightarrow$(4)  Suppose (2).  Since $\rank N \leq 1$ and $N\geq 0$ there exists a $1\times 3n$ vector $\ga$ such that $N=\ga^*\ga$.  Consider the Schur complement identity
\beq \label{SchurCompl}
\bbm A&B\\B^*&D \ebm = \bbm I_2&0\\B^*A\inv &I_{3n}\ebm \bbm A&0\\0&D-B^* A\inv B\ebm \bbm I_2& A\inv B\\ 0&I_{3n} \ebm
\eeq
where $A, D$ are of types $2\times 2$, $3n\times 3n$ respectively.  Choose
\[
A=\bbm -1&0\\0&1\ebm, \quad B=\bbm \ga \\ \ga Z \ebm, \quad D = M-\La^* M\La.
\]
The identity \eqref{SchurCompl} becomes
\begin{align}
\bbm -1 & 0 & \ga \\ 0 & 1 & \ga Z\\  \ga^* &  Z^*\ga^* & M - \Lambda^* M \Lambda \ebm &=~\label{ident}\\
~ & \hspace*{-4cm} \bbm 1 & 0 & 0 \\ 0 & 1 & 0 \\ -\ga^* & Z^*\ga^* & I_{3n} \ebm
\bbm -1 & 0 &  0 \\ 0 & 1 & 0 \\0 & 0 & M - \Lambda^* M  \Lambda  +\ga^*\ga- Z^*\ga^*\ga Z\ebm 
\bbm 1 & 0 & -\ga \\ 0 & 1 & \ga Z\\  0 & 0 & I_{3n}\ebm. \notag
\end{align}
Let
\[
P= - B^*A\inv=\bbm \ga^* & -Z^*\ga^* \ebm \in \C^{3n\times 2}.
\]
Thus equation \eqref{ident} is
\begin{align}
\bbm -1&0&\ga \\ 0&1&\ga Z \\ \ga*&Z^*\ga^*& M-\La^*M\La \ebm &=~\label{ident2}\\
~ & \hspace*{-4cm}
\bbm I_2&0 \\ -P & I_{3n} \ebm
\bbm -1 & 0 & 0 \\ 0 & 1 & 0\\ 0 & 0 & M-\La^*M\La+\ga^*\ga-Z^*\ga^*\ga Z \ebm
\bbm I_2 & -P^* \\ 0 & I_{3n} \ebm. \notag
\end{align}
On pre- and post-multiplying by the inverses of the first and third matrices on the right hand side and using the relation \eqref{Xgeq} we obtain the relation \eqref{1.1bis-LMI}
\begin{align*}
\bbm I_2&0 \\ P & I_{3n} \ebm
\bbm -1&0&\ga \\ 0&1&\ga Z \\ \ga*&Z^*\ga^*& M-\La^*M\La \ebm
\bbm I_2 & P^* \\ 0 & I_{3n} \ebm  &=~\\
~ & \hspace*{-4cm}
\bbm -1 & 0 & 0 \\ 0 & 1 & 0\\ 0 & 0 & M-\La^*M\La+\ga^*\ga-Z^*\ga^*\ga Z \ebm 
	\leq \bbm -1 & 0 & 0 \\ 0 & 1 & 0\\ 0 & 0 & X \ebm.
\end{align*}
Hence (4) holds.  Thus (2)$\Rightarrow$(4); the proof that (3)$\Rightarrow$(5) is almost identical.\\

\noindent (4)$\Rightarrow$(2)  Suppose (4).  The inequality \eqref{1.1bis-LMI} can be written
\[
\bbm I_2 &0 \\ -P& I_{3n} \ebm 
\bbm A&0\\0&X \ebm 
\bbm I_2& -P^*\\ 0& I_{3n} \ebm \geq \bbm A & B \\ B^* & M-\La^* M \La \ebm.
\]
It follows that
\begin{align*}
0\quad &\leq \quad 
\bbm A&-AP^*\\ -PA & PAP^*+X \ebm
 - \bbm A&B\\B^*& M-\La^* M \La \ebm \\
	&=\quad \bbm 0 & -AP^*- B \\ -PA-B^* & PAP^*+X-M+\La^* M\La \ebm.
\end{align*}
Hence $P^* =-A B$ and
\begin{align*}
0 &\leq PAP^*+X-M+\La^* M\La \\
	&=B^*A^3B+X-M+\La^* M\La\\
	&=\bbm \ga^*&Z^*\ga^*\ebm \bbm -1&0\\0&1\ebm \bbm \ga\\\ga Z\ebm+X-M+\La^* M\La\\
	&=-\ga^*\ga+Z^*\ga^*\ga Z+X-M+\La^* M\La
\end{align*}
and so (2) holds with $N=\ga^*\ga$.  Again, (5)$\Rightarrow$(3) is proved in much the same way.
\end{proof}
Relaxation of the condition $\rank N \leq 1$ in (2) of Theorem \ref{NPspectral-LMI} yields a necessary condition for the solvability of a spectral Nevanlinna-Pick problem in the form of the feasibility of an LMI.  The following statement is immediate from Theorem \ref{NPspectral-LMI}.

\begin{corollary}\label{necessary}
In the notation of Theorem {\rm\ref{NPspectral-LMI}}, if the spectral Nevanlinna-Pick problem \eqref {interpMat-LMI}-\eqref{specCond-LMI} is solvable then there exist positive $3n$-square matrices $N$ and $M$ such that the inequality \eqref{Xgeq} holds.
\end{corollary}
In fact the existence of positive $N$ and $M$ such that the inequality \eqref{Xgeq} holds is equivalent to the existence of $\ph\in\schur$ such that
\[
\ph(z_\ell, \la_j) =\Phi(z_\ell, s_j,p_j), \quad \mbox{ for } \ell=1,2,3, \, j=1,\dots,n.
\]
Since $\ph(\cdot,\la)$ need not be linear fractional, we cannot derive an $h\in\hol(\D,\Gamma)$ from $\ph$.

\section{Construction of all interpolating functions}\label{all}
Theorem \ref{NPspectral-criterion} gives us a criterion for the solvability of the interpolation problem
\beq\label{GaInterp}
\textit{find } h\in\hol(\D,\Gamma) \textit{ such that } h(\la_j)=(s_j,p_j) \textit{ for } j=1,\dots,n.
\eeq
The  proof of the theorem contains a description of a process for the derivation of a solution of the problem \eqref{GaInterp} from a feasible pair $(N,M)$ for the LMI  \eqref{Xgeq} or \eqref{inq1.1} with $\rank N \leq 1$.   The process can be summarized as follows.
\begin{center}  \bf Procedure SW \end{center}
Let $z_1, z_2, z_3$ and $\la_j, s_j, p_j$ be as in Theorem \ref{NPspectral-criterion}.  Let $N, M $ be positive $3n$-square matrices such that $\rank N \leq 1$ and the LMI \eqref{Xgeq} holds.
\begin{enumerate}
\item  Choose scalars $\ga_{jk}$ such that $N=\bbm \overline{\ga_{i\ell}}\ga_{jk}\ebm_{i,j=1,\ell,k=1}^{n,3}$.
\item Choose a Hilbert space $\M$ and vectors $v_{jk}\in\M$ such that $M=\bbm\ip{v_{jk}}{v_{i\ell}}_\M\ebm_{i,j=1,\ell,k=1}^{n,3}$.
\item Choose a contraction
\[
\bbm A&B\\C&D \ebm : \C^2\oplus\M \to \C^2\oplus\M
\]
such that
\beq\label{P4.1}
\bbm A&B\\C&D \ebm \bpm 1\\ z_k\ga_{jk} \\ \la_j v_{jk} \epm = \bpm -\Phi(z_k,s_j,p_j) \\ \ga_{jk} \\ v_{jk} \epm
\eeq
for $j=1,\dots,n$ and $k=1,2,3$.
\item Let
\beq\label{thisish}
h(\la) = (\tr,\det)(A+B\la(I-D\la)\inv C)
\eeq
for $\la\in\D$.
\end{enumerate}
Then $h\in\hol(\D,\Gamma)$ and $h(\la_j)=(s_j,p_j)$ for $j=1,\dots,n$.

The purpose of this section is to show that this procedure in principle yields the {\em general} solution of the problem \eqref{GaInterp}, provided that one can find the general feasible pair $(N,M)$ for the relevant LMI with $\rank N \leq 1$.
\begin{proposition}\label{every}
{\em Every} solution of a $\Gamma$-interpolation problem arises by Procedure SW from a solution $(N,M)$ of the corresponding LMI with $\rank N \leq 1$.
\end{proposition}
\begin{proof}
Let $z_k, s_j, p_j$ be as Theorem  \ref{NPspectral-criterion} and let $h\in\hol(\D,\Gamma)$ satisfy $h(\la_j)=(s_j,p_j)$ for $j=1,\dots,n$.  We must produce a pair of positive matrices $(N,M)$ that satisfy the LMI \eqref{Xgeq} such that Procedure SW, when applied to $(N,M)$ with appropriate choices, produces $h$.

By Proposition \ref{realiz-Gamma} there is a unique $F\in \mathcal{S}^{2\times 2}$ such that $h=(\tr F, \det F), \, F_{11} = F_{22}, \, |F_{12}|=|F_{21}|$ a.e. on $\T, \, F_{12}$ is either $0$ or outer and $F_{12}(0) \geq 0$.  Then  $F_{11}=F_{22}=\half s$, and  Proposition \ref{realiz-Gamma} asserts further that if 
\begin{align*}
\ga(z,\la) &= \frac{F_{12}(\la)}{1-\half s(\la) z}, \\
\eta(z,\la) &= \bbm 1 \\ \ga(z,\la) \ebm \quad \mbox{ and } \\
J&= \bbm 0&1 \\ 1& 0 \ebm
\end{align*}
then
\[
1-\overline{\Phi(w,h(\mu))}\Phi(z,h(\la)) = (1-\bar w z)\overline{\ga(w,\mu)}\ga(z,\la) + \eta(w,\mu)^* J\left(I-F(\mu)^*F(\la)\right)J\eta(z,\la)
\]
for all $z,w,\la,\mu\in \D$.  Since $F\in\mathcal{S}^{2\times 2}$, the map
\[
(\la,\mu) \mapsto  J\frac{I-F(\mu)^*F(\la)}{1-\bar\mu \la} J
\]
is a positive $2\times 2$ kernel on $\D$, and so there is a Hilbert space $\h$ and an analytic map
$U: \D\to \mathcal{L}(\C^2, \h)$ such that
\[
J\frac{I-F(\mu)^*F(\la)}{1-\bar\mu \la} J = U(\mu)^* U(\la)
\]
for all $\la,\mu\in\D$.  Then
\beq\label{P7.1}
1-\overline{\Phi(w,h(\mu))}\Phi(z,h(\la)) = (1-\bar w z)\overline{\ga(w,\mu)}\ga(z,\la) + (1-\bar\mu\la)\eta(w,\mu)^*  U(\mu)^* U(\la)\eta(z,\la).
\eeq
In particular, when $w=z_\ell, \ \mu=\la_i, \ z=z_k, \ \la=\la_j$, 
\begin{align*}
1-\overline{\Phi(z_\ell,h(s_i,p_i)}\Phi(z_k,s_j,p_j) &= (1-\bar z_\ell z_k)\overline{\ga(z_\ell,\la_i)}\ga(_kz,\la_j)\\
	& \hspace*{-3cm} + (1-\bar\la_i\la_j) \ip{U(\la_j)\eta(z_k,\la_j)}{U(\la_i)\eta(z_\ell,\la_i)}_\h
\end{align*}
for $i,j=1,\dots,n, \, \ell,k=1,2,3$.  Thus the $3n$-square matrices
\begin{align*}
N&= \bbm \overline{\ga(z_\ell,\la_i)}\ga(_kz,\la_j) \ebm \\
M&= \bbm \ip{U(\la_j)\eta(z_k,\la_j)}{U(\la_i)\eta(z_\ell,\la_i)}_\h \ebm
\end{align*}
satisfy the LMI \eqref{Xgeq} (and even the matrix equation \eqref{1.1}), and $\rank N \leq 1$.  We may therefore apply Procedure SW to $(N,M)$.   In steps (1) and (2) choose
\[
\ga_{jk} = \ga(z_k,\la_j), \quad \M=\h,\quad v_{jk} = U(\la_j)\eta(z_k,\la_j).
\]
By virtue of the relation \eqref{P7.1} the Gramian of the vectors
\beq\label{spanning}
\bpm 1 \\ z\ga(z,\la)\\ \la U(\la) \eta(z,\la) \epm \in \C^2\oplus \h, \qquad z,\la \in \D,
\eeq
is equal to the Gramian of the vectors
\[
\bpm -\Phi(z,h(\la)) \\ \ga(z,\la) \\ U(\la) \eta(z,\la) \epm  \in \C^2\oplus \h, \qquad z,\la \in \D.
\]
Hence there exists an isometry $L_0$ on the subspace of $\C^2\oplus\h$ spanned by the vectors \eqref{spanning} such that
\beq\label{P9.1}
L_0 \bpm 1 \\ z\ga(z,\la)\\ \la U(\la) \eta(z,\la) \epm = \bpm -\Phi(z,h(\la)) \\ \ga(z,\la) \\ U(\la) \eta(z,\la) \epm
\eeq for all $z,\la \in\D$.  Let 
\[
L=\bbm A&B\\C&D \ebm \in \mathcal{L}(\C^2\oplus\h)
\]
be any contractive extension of $L_0$.  On specialising equation \eqref{P9.1} to $z_k$ and $\la_j$ one obtains the relation \eqref{P4.1} in step 3 of Procedure SW.  One may therefore use $L$ in step 4, and so obtain a function $\tilde h\in\hol(\D,\Gamma)$ that satisfies $\tilde h(\la_j)=(s_j,p_j)$.

We claim that $\tilde h = h$.  By equation \eqref{P9.1},
\begin{align*}
\bpm -\Phi(z(h(\la)) \\ \ga(z,\la) \epm &= A \bpm 1\\z\ga(z,\la) \epm + B\la U(\la)\eta(z,\la) \\
U(\la)\eta(z,\la) &= C\bpm 1 \\ z\ga(z,\la) \epm + D\la U(\la) \eta(z,\la)
\end{align*}
 and so, by elimination of $\eta(z,\la)$,
\begin{align*}
\bpm -\Phi(z,h(\la) \\ \ga(z,\la) \epm &= \left( A+B\la(I-D\la)\inv C\right) \eta(z,\la) \\
	&=\Psi(\la) \bpm 1 \\ z\ga(z,\la) \epm
\end{align*}
for all $z,\la\in\D$.  Now eliminate $\ga(z,\la)$ to obtain
\[
-\Phi(z,h(\la)) = \Psi_{11}(\la) + \frac{\Psi_{12}\Psi_{21}(\la)z}{1-\Psi_{22}(\la)z}
\]
for all $z,\la\in\D$.
Since
\[
-\Phi(z,h(\la)) = \half s(\la) + \frac{\left(\tfrac 14 s(\la)^2-p(\la)\right) z}{1-\half s(\la) z}
\]
it follows that
\[
\Psi_{11}(\la)= \half s(\la) = \Psi_{22}(\la)
\]
and
\[
\Psi_{12}(\la)\Psi_{21}(\la) = \tfrac 14 s(\la)^2-p(\la).
\]
Hence
\[
\tr \Psi = s, \qquad \det \Psi = p
\]
and therefore $\tilde h = h$ as required.
\end{proof}

\section{Implementation of the solution procedure}\label{implement}

We conclude with some remarks on the practical feasibility of our results for the numerical solution of a spectral Nevanlinna-Pick problem.
Let interpolation points $\la_1,\dots, \la_n$ and target matrices $W_1,\dots,W_n$ of type $2\times 2$ be given.  If any of the $W_j$ are scalar matrices then the corresponding interpolation conditions can be removed by the standard process of Schur reduction, and so we may suppose that all the $W_j$ are nonscalar.  Alternatively, if some $W_j$ are scalar, one may still reduce to an interpolation problem for $\hol(\D,\Gamma)$, but with interpolation conditions on derivatives \cite{NiPfTh}; one could then try to prove analogs of the present results for this wider class of interpolation problems (this should not be difficult).

Supposing, then, that the $W_j$ are nonscalar, let $s_j=\tr W_j, \, p_j=\det W_j$.  As recalled in Proposition \ref {reducetoGamma}, the problem reduces to the solution of the interpolation problem \eqref{GaInterp}.  Choose $z_1,z_2, z_3$ of modulus at most $1$ such that no $(s_j,p_j)$ is $(2z_k,z_k^2)$ for any $k$ (if $r(W_j) < 1$ for each $j$ then one can make the natural choice of $-1,0,1$ for the $z_k$).

To determine with the aid of Theorem \ref{NPspectral-main} whether the problem \eqref{GaInterp} with these data is solvable  we may test the  criterion (3) of the theorem.  That is, we must ascertain whether there exist positive matrices $N$ of rank $1$ and $M$  satisfying the LMI \eqref{object} in condition (2).  Existing software packages can reliably determine whether such an LMI is feasible, but we do not know an effective way to test whether there is a feasible pair such that $\rank N \leq 1$.  The following refinement of Theorem \ref{NPspectral-criterion} at least shows that a search over a {\em compact} set of pairs $(N,M)$ suffices.

\begin{proposition}\label{boundsNM}
Let $\la_j, s_j, p_j$ and $z_k$ be as in Theorem {\rm \ref{NPspectral-criterion}}. The $\Gamma$-interpolation problem  
\beq\label{theproblem}
\la_j\in\D  \mapsto (s_j,p_j)\in\Gamma, \qquad j=1,\dots, n,
\eeq
is solvable if and only if 
there exist positive $3n$-square matrices  $N= [N_{i\ell,jk}]_{i,j=1,\, \ell,k=1}^{n,3}$ of rank $1$ and $M=[M_{i\ell,jk}]_{i,j=1,\, \ell,k=1}^{n,3}  $  that satisfy the LMI \eqref{Xgeq} and 
\begin{eqnarray} 
\label{M-bounded}
|M_{i\ell,jk}| & \le & \frac{2}{|1 - \bar{\la_i} \la_j|} \sqrt{1 + \frac{1}{(1 - 
\half |s_j|)^2}}\sqrt{1 + \frac{1}{(1 - 
\half |s_j|)^2}},
\end{eqnarray}
and 
\begin{eqnarray} 
\label{N-bounded}
|N_{i\ell,jk}| & \le & \frac{1}{(1 - 
\half |s_i|)(1 - \half |s_j|)}.
\end{eqnarray}
\end{proposition}
\begin{proof}
Sufficiency is contained in Theorem \ref{NPspectral-criterion}, (3)$\Rightarrow$(1).  To prove necessity, suppose that the interpolation problem is solvable.  In the proof of Theorem \ref{NPspectral-criterion} (1)$\Rightarrow$(3) it was shown that the LMI \eqref {inq1.1}  (or equivalently, \eqref{Xgeq}) holds when
\begin{eqnarray} \label{N-bounds}
N = [N_{i\ell,jk}]_{i,j=1,\, \ell,k=1}^{n,3} = \left[\overline{\gamma(\la_i,z_\ell)}\gamma(\la_j,z_k)\right]_{i,j=1,\, \ell,k=1}^{n,3} 
\end{eqnarray}  
of rank $1$ and 
\begin{eqnarray} 
\label{M-bounds}
M=[M_{i\ell,jk}]_{i,j=1,\, \ell,k=1}^{n,3} = \left[
\eta(\la_i, z_\ell)^* \frac{ I - F(\la_i)^* F(\la_j)}{1- \bar{\la_i} \la_j}\eta(\la_j, z_k)\right]_{i,j=1,\, \ell,k=1}^{n,3}  
\end{eqnarray} 
where each $F(\la_j)$ is a contraction,
\begin{eqnarray} 
\gamma(\la_j,z_k)&=&(1 - \half s_j z_k)^{-1} f_1(\la_j), \;\;  \nonumber\\
\eta(\la_j, z_k)&=& \bbm \gamma(\la_j,z_k)z_k\\ 1 \ebm
\end{eqnarray} 
and the function $f_1$ is in the Schur class.
Thus, for $j =1, \dots, n$ and $k=1,2,3$,
\beq 
\label{gamma-bounded}
 |\gamma(\la_j, z_k)|\leq \frac{1}{1 - \half |s_j|},
\eeq
from which the estimate \eqref{N-bounded} follows.  Moreover
\[
\|\eta(\la_j, z_k)\|^2 =  \left\| \bbm \gamma(\la_j,z_k)z_k\\ 1 \ebm \right\|^2 \le  1 + \frac{1}{(1 - \half |s_j|)^2}
\]
and therefore
\[
|M_{i\ell,jk}| \leq \frac{\|\eta(\la_i, z_\ell)\|\, \|\eta(\la_j, z_k)\|}{|1-\bar\la_i \la_j|} \|I-F(\la_i)^*F(\la_j)\|,
\]
from which the bound \eqref{M-bounded} follows.
\end{proof}

One approach to the finding of a suitable pair $(N,M)$ would be to use the alternative formulation (2)
 in Theorem \ref{NPspectral-criterion}: to minimize the quadratic function $f(N, M)=(\tr N)^2-\tr(N^2)$ over the feasible region $R$.   By compactness $f$ attains its minimum on the set of positive pairs $(N,M)$ that satisfy the LMI \eqref{Xgeq} and the bounds \eqref{N-bounded} and \eqref{M-bounded}, provided that this set is nonempty.  Proposition \ref{boundsNM} asserts that the $\Gamma$-interpolation problem \eqref{theproblem} is solvable if and only if this minimum is zero.

Since $f$ is positive homogeneous of degree $2$ its local minima over $R$ all lie on the topological boundary of $R$, and the gradient of $f$ is linear in $N$.  However, $R$ is a subset of real Euclidean space of $18n^2$ variables and the boundary of $R$ has a complicated structure.

Once a feasible pair $(N,M)$  for the LMI \eqref{Xgeq} with $\rank N \leq 1$ is found, it is a matter of straightforward linear algebra to apply Procedure SW in order to calculate a solution $h$ of the interpolation problem \eqref{theproblem}.
It is routine to find $3n$ scalars $\ga_{jk}$ such that $N_{i\ell,jk} = \bar \ga_{i\ell}\ga_{jk}$ for all $i,\ell,j,k$.  Likewise, by Cholesky factorization, one can find $3n$ vectors $v_{jk}$ in some Hilbert space $\M$ such that $M_{i\ell,jk} =\ip{v_{jk} }{v_{i\ell}}$.  Because the LMI \eqref{Xgeq} holds, there is an isometric operator matrix $\bbm A&B\\C&D \ebm$ that satisfies the relation \eqref{P4.1}; one may then define an interpolating $\Gamma$-inner function $h$ by equation \eqref{thisish}.   In principle it is  simple linear algebra to find $A,B,C$ and $D$ when $\M$ is chosen to be finite-dimensional.

As Proposition \ref{every} shows, the above procedure when applied to the general feasible pair $(N,M)$ for the LMI \eqref{Xgeq} with $\rank N \leq 1$  yields all possible interpolating functions.  For numerical implementation one would naturally take $\M$ finite-dimensional, and then the resulting function $h$ will be rational.  A slight modification of Proposition \ref{every} shows that all {\em rational} interpolating functions are obtainable by this procedure. 

An important question is whether the results of this paper furnish an improvement (for the special case of the spectral Nevanlinna-Pick problem) on existing `$D$-$K$ iteration' methods, as for instance in the Matlab mu-analysis toolbox \cite{mathworks}.  We leave this question for future exploration.

\section*{Acknowledgements} 

The first author was partially supported by National Science Foundation Grant on  Extending Hilbert Space Operators DMS 1068830.
The third author was partially supported by the UK Engineering and Physical Sciences Research Council grants EP/J004545/1 and EP/K50340X/1. We are grateful to Professors R. Bitmead and R. Skelton of the University of California at San Diego for helpful suggestions.

JIM AGLER, Department of Mathematics, University of California at San Diego, San Diego, CA \textup{92103}, USA\\

ZINAIDA A. LYKOVA,
School of Mathematics and Statistics, Newcastle University, Newcastle upon Tyne
 NE\textup{1} \textup{7}RU, U.K.~~\\
e-mail\textup{: \texttt{Zinaida.Lykova@ncl.ac.uk}}\\

N. J. YOUNG, School of Mathematics, Leeds University, Leeds LS2 9JT, U.K.~~
and School of Mathematics and Statistics, Newcastle University, Newcastle upon Tyne
 NE\textup{1} \textup{7}RU, U.K.~~\\
e-mail\textup{: \texttt{N.J.Young@leeds.ac.uk}}\\
\end{document}